\documentclass[showkeys,amsmath,amssymb]{revtex4-2}

\usepackage{graphicx}% Include figure files
\newcommand{\p}{\partial}
\newcommand{\be}{\begin{equation}}
\newcommand{\ee}{\end{equation}}

\begin{document}

%\preprint{APS/123-QED}

\title{Chimeras with uniformly distributed heterogeneity: two coupled populations}% Force line breaks with \\

\author{Carlo R. Laing}
 \email{c.r.laing@massey.ac.nz}
\affiliation{School of Natural and Computational Sciences, 
Massey University,
Private Bag 102-904 
North Shore Mail Centre, 
Auckland,
New Zealand
}%

\date{\today}% It is always \today, today,

\begin{abstract}
Chimeras occur in networks of two coupled populations of oscillators when the oscillators
in one population synchronise while those in the other are asynchronous. We consider 
chimeras of this form in networks
of planar oscillators for which one parameter associated with the dynamics of an oscillator
is randomly chosen from a uniform distribution. A generalisation of the approach
in [C.R. Laing, Physical Review E, {\bf 100}, 042211, 2019], which dealt with identical
oscillators, is used to investigate the existence and stability of chimeras for these
heterogeneous networks in the limit of an infinite number of oscillators. In all cases,
making the oscillators more heterogeneous destroys the stable chimera in a saddle-node
bifurcation. The results help us understand the robustness of chimeras in networks of general
oscillators to heterogeneity.

\end{abstract}

%Use showkeys class option if keyword
\maketitle

\section{Introduction}
Chimera states occur in networks of coupled oscillators and are characterised by coexisting
groups of synchronised and asynchronous oscillators~\cite{panabr15,ome18}. They have been
observed in one-dimensional~\cite{abrstr06,omezak15} and 
two-dimensional~\cite{lai17,panabr15a} domains with nonlocal coupling, and also networks
formed from two populations with strong coupling within a population and weaker coupling
between them~\cite{abrmir08,pikros08,panabr16}. A variety of oscillator types have been
considered, with the most common being a phase oscillator~\cite{kurbat02}, but others include
Stuart-Landau oscillators~\cite{lai10}, van der Pol oscillators~\cite{omezak15,uloome16},
oscillators with inertia~\cite{boukan14,olm15} and neural 
models~\cite{olmpol11,omeome13,ratpyr17}.

Many investigations of chimeras report only the results of numerically solving a finite
number of ordinary differential equations (ODEs) which describe the networks' behaviour.
Such simulations are for only a finite time, so the results seen may actually be part of a long
transient~\cite{Zakkap16}. 
With a finite network there is the issue of finite size effects, such as positive
Lyapunov exponents which tend to zero as the network size is increased~\cite{ome18}
or chimeras' finite lifetimes~\cite{wolome11}. Perhaps most significantly, such simulations
cannot detect unstable states so it is often not clear what happens to a stable chimera
as a parameter is varied, other than it no longer existing.

Early results on the existence of chimeras used a self-consistency 
approach~\cite{lai10,kurbat02,shikur04,abrstr04} but this does not provide information on the stability of solutions. A great deal of progress has been made using the Ott/Antonsen
ansatz~\cite{ottant08,ottant09}, since it gives evolution equations for quantities
of interest, but its use is restricted to networks of phase oscillators
coupled through sinusoidal functions of phase differences~\cite{abrmir08,lai09,ome18}.
Laing~\cite{lai10} used self-consistency to investigate the existence of chimeras in networks
of two populations of Stuart-Landau oscillators, each oscillator being 
described by a complex variable. This was later generalised~\cite{lai19} using techniques
from~\cite{clupol18} to determine the stability of these chimeras, and chimeras in networks
of three more types of oscillators (Kuramoto with inertia, FitzHugh-Nagumo oscillators,
delayed Stuart-Landau oscillators) were studied.

The approach in~\cite{lai19} was to recognise that the incoherent oscillators in one
population lie on a curve $\mathcal{C}$ in the phase plane while those in the synchronous
population can be described by a pair of real variables, since all of these oscillators are
identical and undergo the same dynamics. In the limit of an infinite number
of oscillators in each population the curve $\mathcal{C}$ is described by its
shape (distance from the origin in polar coordinates) and the density of oscillators on it,
and partial differential equations (PDEs) governing the evolution of these functions
can be derived~\cite{clupol18}. The full network is then described by a pair of PDEs and a 
pair of ODEs, coupled by an integral.

The analysis in~\cite{lai19} assumed {\em identical} oscillators, but we do not expect this
to be the case in any experimental situation~\cite{tinnko12,marthu13,totrod18} and it is
known that networks of identical oscillators may have qualitatively different dynamics
from those of nonidentical oscillators~\cite{watstr94}. In this paper we extend the
results in~\cite{lai19} to the case of nonidentical oscillators. Specifically, we assume
that for each oscillator, one parameter associated with its dynamics
is randomly chosen from a uniform distribution. A uniform distribution is zero outside some
range and this means that for a narrow distribution, the types of chimeras observed
in~\cite{lai19} persist and can be described by a generalisation of the techniques
developed in that paper.

Various distributions of intrinsic frequencies in networks of all-to-all 
coupled oscillators have been
considered, e.g.~Lorentzian~\cite{ottant08}, 
bimodal~\cite{marbar08}, Gaussian~\cite{strmir91,hanfor18}, beta~\cite{daeds20}
and uniform~\cite{piedes18,paz05,bairos10,ottstr16,eydwol17}. 
There are significant differences in the transition
to synchrony as coupling strength is increased between distributions with compact support
and those whose support is unbounded. In the former case one normally observes a first-order
transition, whereas in the latter it is second-order. Also, for an infinite network
full synchrony --- in which all oscillators are phase locked --- can only occur when the
frequency distribution has compact support~\cite{erm85}. We observe and exploit this phenomenon
to analyse the networks studied in this paper.

Previous relevant work includes~\cite{rybvad19}, which considers a network formed from coupled
ring subnetworks of logistic maps in which a number of parameters are randomly chosen from uniform
distributions. The authors investigate the effects of varying the widths of these distributions
on the number of subnetworks which fully synchronise.

We consider networks formed from two populations of oscillators.
In Sec.~\ref{sec:kur} we consider Kuramoto-type phase oscillators and in Sec.~\ref{sec:SL}
we revist the Stuart-Landau oscillators studied in~\cite{lai10}. Sec~\ref{sec:pend}
considers Kuramoto oscillators with inertia, also studied in~\cite{lai10}.
We study van der Pol oscillators in Sec.~\ref{sec:vdp} and conclude in Sec.~\ref{sec:disc}.

\section{Kuramoto phase oscillators}
\label{sec:kur}
We first consider two populations of phase oscillators coupled through a sinusoidal
function of phase differences. Networks of this form have
been studied previously~\cite{abrmir08,lai09a,panabr16,pikros08,monkur04}.
We first consider heterogeneity in intrinsic frequencies, then in the strength of
coupling between populations. 

\subsection{Distributed frequencies}
\label{sec:varomkur}

Consider two populations of $N$ phase oscillators each governed by
\be
   \frac{d\theta_j}{dt}  = \omega_j+\frac{\mu}{N}\sum_{k=1}^N\sin{(\theta_k-\theta_j-\alpha)}+\frac{\nu}{N}\sum_{k=1}^N\sin{(\theta_{N+k}-\theta_j-\alpha)} \label{eq:dth1}
\ee
for $j=1,2\dots N$
and
\be
   \frac{d\theta_j}{dt}  = \omega_j+\frac{\mu}{N}\sum_{k=1}^N\sin{(\theta_{N+k}-\theta_j-\alpha)}+\frac{\nu}{N}\sum_{k=1}^N\sin{(\theta_{k}-\theta_j-\alpha)} \label{eq:dth2}
\ee
for $j=N+1,N+2,\dots 2N$. $\mu$ is the strength of coupling within a population and
$\nu$ is the strength between populations. 
For identical $\omega_j$ this system reduces to the system studied
in~\cite{abrmir08,panabr16,pikros08} while if they are chosen from a Lorentzian 
distribution it is the same as 
in~\cite{lai09a}.
Instead, here for each population the $\omega_j$ are randomly chosen from the
uniform distribution $p(\omega)$ which is non-zero only on the interval $B$.

\begin{figure}
\begin{center}
\includegraphics[width=13cm]{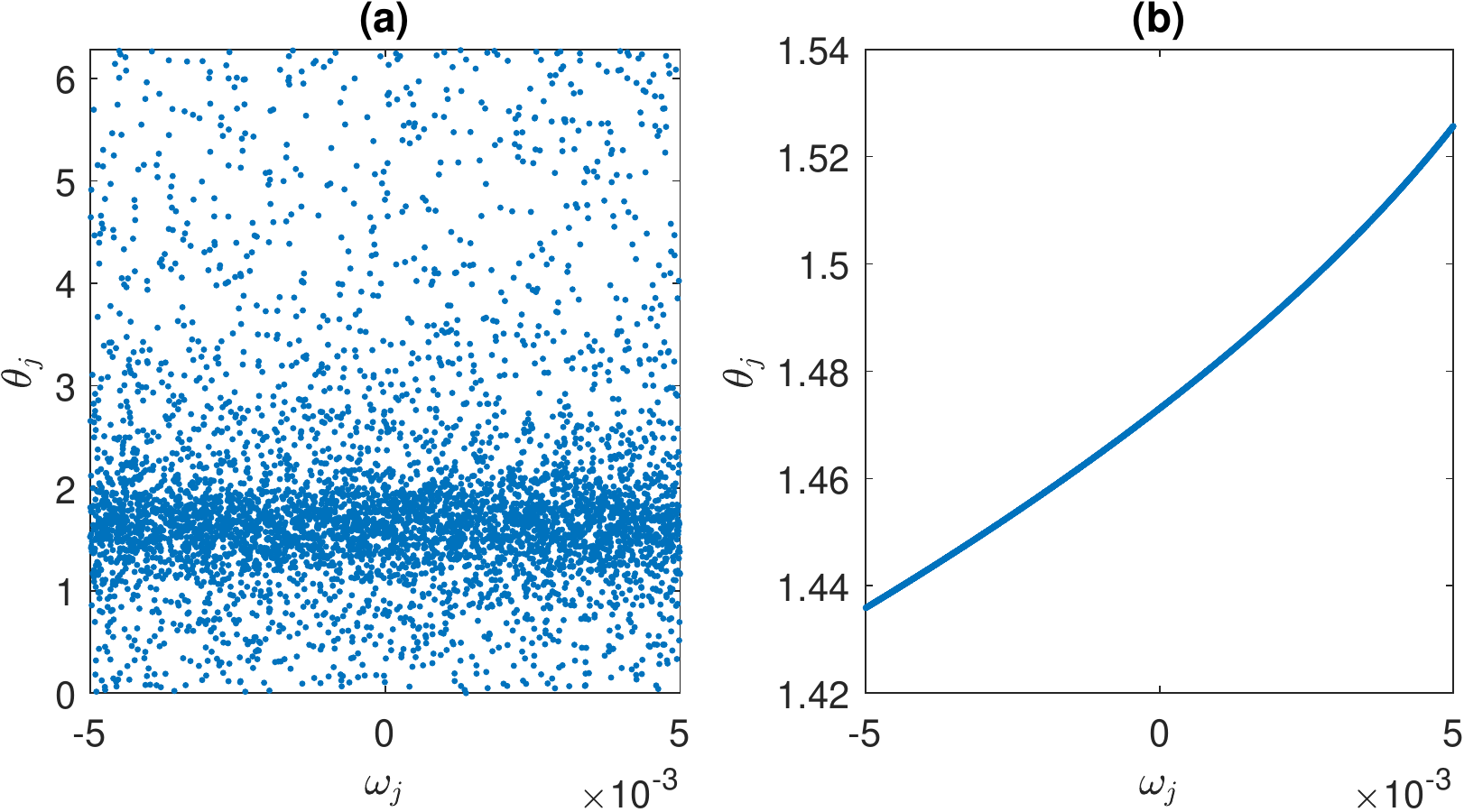}
\caption{A snapshot of a chimera state solution of~\eqref{eq:dth1}-\eqref{eq:dth2}.
(a): population 1; (b): population 2. Note the different vertical axes.
Parameters: $N=5000,\mu=0.6,\nu=0.4,\alpha=\pi/2-0.08,\Delta\omega=0.005$.}
\label{fig:snapkur}
\end{center}
\end{figure}

An example of a chimera state for~\eqref{eq:dth1}-\eqref{eq:dth2}
is shown in Fig.~\ref{fig:snapkur} where $p(\omega)$
is uniform on $[-\Delta\omega,\Delta\omega]$. We see that population 1 is incoherent,
with no apparent dependence of $\theta_j$ on $\omega_j$, whereas population 2 is synchronised
(although not phase synchronised) with a clear dependence of $\theta_j$ on $\omega_j$.
The average frequencies of oscillators in the two populations are different, as required
for a chimera state. [This state is close to the one which occurs for identical oscillators,
so we also refer to it (and many states studied below) as a ``chimera''.]
We now proceed to analyse this state, in terms of both existence and stability.

We assume that population 2 is locked, 
and write $\theta_{N+j}=\phi_j$ for $j=1,2,\dots N$. Thus, using trigonometric
identities, for population 2
\begin{align}
   \frac{d\phi_j}{dt}  & = \omega_{N+j}+\frac{\mu}{N}\sum_{k=1}^N\sin{(\phi_{k}-\phi_j-\alpha)}+\frac{\nu}{N}\sum_{k=1}^N\sin{(\theta_{k}-\phi_j-\alpha)} \nonumber \\
 % &  =\omega_{N+j}+\mu\left[\cos{(\phi_j+\alpha)}\hat{S}-\sin{(\phi_j+\alpha)}\hat{C}\right]
% +\nu\left[\cos{(\phi_j+\alpha)}S-\sin{(\phi_j+\alpha)}C\right] \\
  & = \omega_{N+j}+(\mu\hat{S}+\nu S)\cos{(\phi_j+\alpha)}-(\mu\hat{C}+\nu C)\sin{(\phi_j+\alpha)} \label{eq:dphidtB}
\end{align}
for $j=1,2,\dots N$ where
\be
   \hat{S}\equiv \frac{1}{N}\sum_{k=1}^N\sin{\phi_k}; \qquad \hat{C}\equiv \frac{1}{N}\sum_{k=1}^N\cos{\phi_k}; \qquad S\equiv \frac{1}{N}\sum_{k=1}^N\sin{\theta_k}; \qquad C\equiv \frac{1}{N}\sum_{k=1}^N\cos{\theta_k}
\ee
For population 1 we have
\begin{align}
   \frac{d\theta_j}{dt} & =\omega_j+\mbox{Im}\left[e^{-i(\theta_j+\alpha)}\left(\frac{\mu}{N}\sum_{k=1}^N e^{i\theta_k}+\frac{\nu}{N}\sum_{k=1}^N e^{i\phi_k}\right)\right] \nonumber \\
 %  & = \omega_j+\mbox{Im}\left[e^{-i(\theta_j+\alpha)}\left(\mu(C+iS)+\nu(\hat{C}+i\hat{S})\right)\right] \\
   & = \omega_j+\frac{1}{2i}\left[e^{-i(\theta_j+\alpha)}Z-e^{i(\theta_j+\alpha)}\bar{Z}\right]
\end{align}
where
\be
   Z\equiv \mu(C+iS)+\nu(\hat{C}+i\hat{S}) \label{eq:Z}
\ee
and overline indicates complex conjugate.
We now take the continuum limit $N\to\infty$. 
For population 2, instead of individual oscillators with phases
$\phi_j$ and frequencies $\omega_j$, $\omega$ is now a continuous parameter and we have the
function $\phi(\omega,t)$, defined for $\omega\in B$.  It satisfies the continuum version 
of~\eqref{eq:dphidtB}:
\be
   \frac{\partial \phi(\omega,t)}{\partial t}=\omega+(\mu\hat{S}+\nu S)\cos{(\phi+\alpha)}-(\mu\hat{C}+\nu C)\sin{(\phi+\alpha)}
\ee
Fig.~\ref{fig:snapkur}(b) can be regarded as showing $\phi(\omega,t)$ for the discrete
values of $\omega$ used in the simulation, for a particular value of $t$.
 
In this limit we have
\be
   \hat{S}=\int_B \sin{[\phi(\omega,t)]}p(\omega)\ d\omega; \qquad \hat{C}=\int_B \cos{[\phi(\omega,t)]}p(\omega)\ d\omega \label{eq:hatS}
\ee

Population 1 is described by the probability density function $f(\theta,\omega,t)$  which
satisfies the continuity equation
\be
   \frac{\partial f}{\partial t}+\frac{\partial}{\partial\theta}(fv)=0
\ee
where
\be
   v(\theta,\omega,t)=\omega+\frac{1}{2i}\left[e^{-i(\theta+\alpha)}Z-e^{i(\theta+\alpha)}\bar{Z}\right]
\ee

We can apply Ott/Antonsen ansatz~\cite{ottant08,ottant09} and write
\be
   f(\theta,\omega,t)=\frac{p(\omega)}{2\pi}\left[1+\sum_{n=1}^\infty a(\omega,t)^ne^{-in\theta}+\mbox{c.c.}\right]
\ee
where ``c.c.'' is the complex conjugate of the previous term.
Substituting this ansatz into the continuity equation above gives the evolution equation for 
$a(\omega,t)$~\cite{ottant08}:
\be
   \frac{\partial a(\omega,t)}{\partial t}= %i\left(\frac{e^{-i\alpha}Z}{2i}+\omega a-\frac{e^{i\alpha}\bar{Z}}{2i}a^2\right)=
    i\omega a+\frac{1}{2}\left(e^{-i\alpha}Z-e^{i\alpha}\bar{Z}a^2\right)
\ee
%where $a(\omega,t)$ is the expected value of $e^{i\theta}$ for the oscillators with frequency
%$\omega$. 
Lastly, 
\be
   C+iS=\int_B\int_0^{2\pi} f(\theta,\omega,t)e^{i\theta}\ d\theta\ d\omega=\int_B a(\omega,t)p(\omega)\ d\omega \label{eq:CiS}
\ee
%We can thus write
%\be
%   Z=\int_B \left[\mu a(\omega,t)+\nu e^{i\phi(\omega,t)}\right] p(\omega)\ d\omega
%\ee
We move to a rotating coordinate frame rotating with  angular speed $\Omega$ in which
both $a$ and $\phi$ are constant. Note that the phases of oscillators in population 1
are not constant in this frame, even though their density is.
Thus we are interested in fixed points of
\begin{align}
    \frac{\partial a(\omega,t)}{\partial t} & =i(\omega-\Omega)a+\frac{1}{2}\left(e^{-i\alpha}Z-e^{i\alpha}\bar{Z}a^2\right) \label{eq:dadt} \\
   \frac{\partial \phi(\omega,t)}{\partial t} & =\omega-\Omega+(\mu\hat{S}+\nu S)\cos{(\phi+\alpha)}-(\mu\hat{C}+\nu C)\sin{(\phi+\alpha)} \label{eq:dphidt}
\end{align}
This is a pair of PDEs, one for the complex quantity $a$ and the other for the angle $\phi$,
coupled through the integrals~\eqref{eq:hatS} and~\eqref{eq:CiS}. The physical interpretation of 
$\phi$ is clear, and for fixed $\omega$, the angular dependence of $f(\theta,\omega,t)$ is a 
Poisson kernel with centre given by the argument of $a(\omega,t)$ and its ``sharpness''
determined by the magnitude of $a(\omega,t)$~\cite{lai09}. Eqns.~\eqref{eq:dadt}-\eqref{eq:dphidt}
can be thought of as a generalisation of eqns.~(11) in~\cite{abrmir08} to the case of
nonidentical oscillators.

Due to the invariance under a global phase shift there is a continuum of fixed points
of~\eqref{eq:dadt}-\eqref{eq:dphidt}, each a
shift of one another. Thus we append a ``pinning'' condition; in this case, $\hat{S}=0$.
This additional equation allows us to find all the unknowns, $a,\phi$ and $\Omega$. 
We choose $p(\omega)$ to be uniform on $[-\Delta\omega,\Delta\omega]$ and use Gauss-Legendre
quadrature with 50 points to approximate the integrals.
Thus the domain $[-\Delta\omega,\Delta\omega]$ is discretised using the points
$\omega_i=\Delta\omega x_i$ for $i=1,2,\dots 50$
where the $x_i$ are the roots of $P_{50}(x)$, the  Legendre polynomial of order 50.
The integrals over $\omega$ in~\eqref{eq:hatS} and~\eqref{eq:CiS} 
are thus approximated by weighted sums.

\subsubsection{Results}
\label{sec:varyomkur}
\begin{figure}
\begin{center}
\includegraphics[width=13cm]{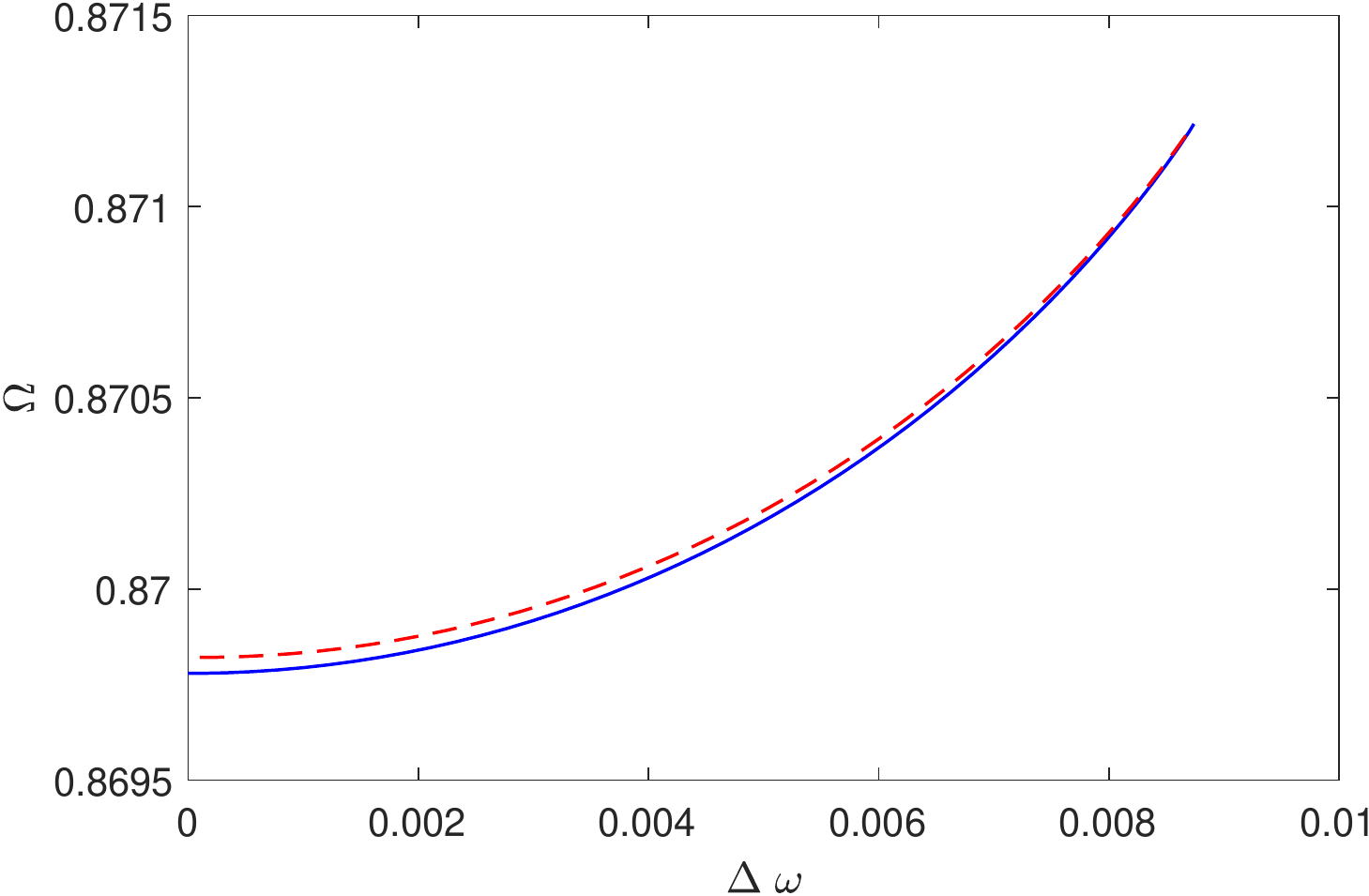}
\caption{$\Omega$ for 
fixed points of~\eqref{eq:dadt}-\eqref{eq:dphidt} describing chimera states. 
Solid: ``stable''. Dashed: unstable.
Parameters: $\mu=0.6,\nu=0.4,\alpha=\pi/2-0.08$.}
\label{fig:vardelomkur}
\end{center}
\end{figure}

The most obvious question is: what is the influence of having distributed values of $\omega$
on the existence and stability of chimeras?
Using pseudo-arclength continuation~\cite{lai14,gov00} and 
varying $\Delta \omega$ we obtain Fig.~\ref{fig:vardelomkur}. The stable chimera 
that exists for identical oscillators is
destroyed in a saddle-node bifurcation as $\Delta \omega$ is increased, i.e.~the oscillators
are made more heterogeneous. This is in contrast to the situation when the $\omega_j$
are chosen from a Lorentzian distribution, where increasing the level of heterogeneity
causes the distribution of phases in the chimera to become more peaked and the distribution
in the synchronous group to be come less peaked until both states meet in a pitchfork
bifurcation and merge to form a state in which the two populations cannot be 
distinguished~\cite{lai09a}.

\begin{figure}
\begin{center}
\includegraphics[width=13cm]{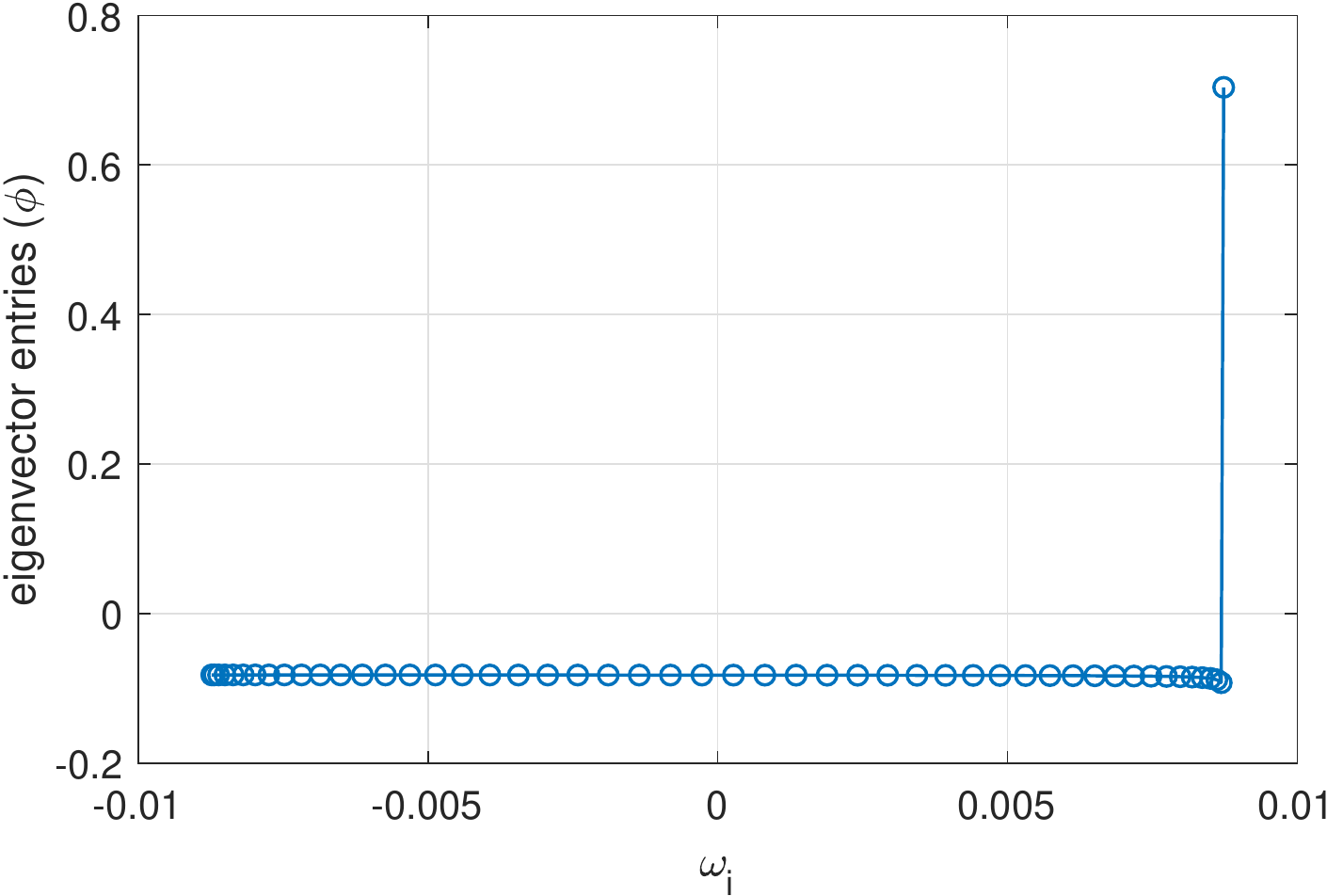}
\caption{Eigenvector localisation.
$\phi$ component of the eigenvector corresponding to the zero eigenvalue at the 
saddle-node bifurcation shown in Fig.~\ref{fig:vardelomkur}.
Parameters: $\mu=0.6,\nu=0.4,\alpha=\pi/2-0.08,
\Delta\omega=0.0087392$.}
\label{fig:eigvec}
\end{center}
\end{figure}

 The $\phi$ component of the eigenvector corresponding to the zero eigenvalue at the 
saddle-node bifurcation is shown in Fig.~\ref{fig:eigvec} and it is clear that this is localised
at the highest $\omega_i$, i.e.~it is the oscillator with the largest intrinsic frequency
which ``unlocks'' first as $\Delta\omega$ is increased, leading to 
``phase walkthrough''~\cite{ermrin84}.

\begin{figure}
\begin{center}
\includegraphics[width=13cm]{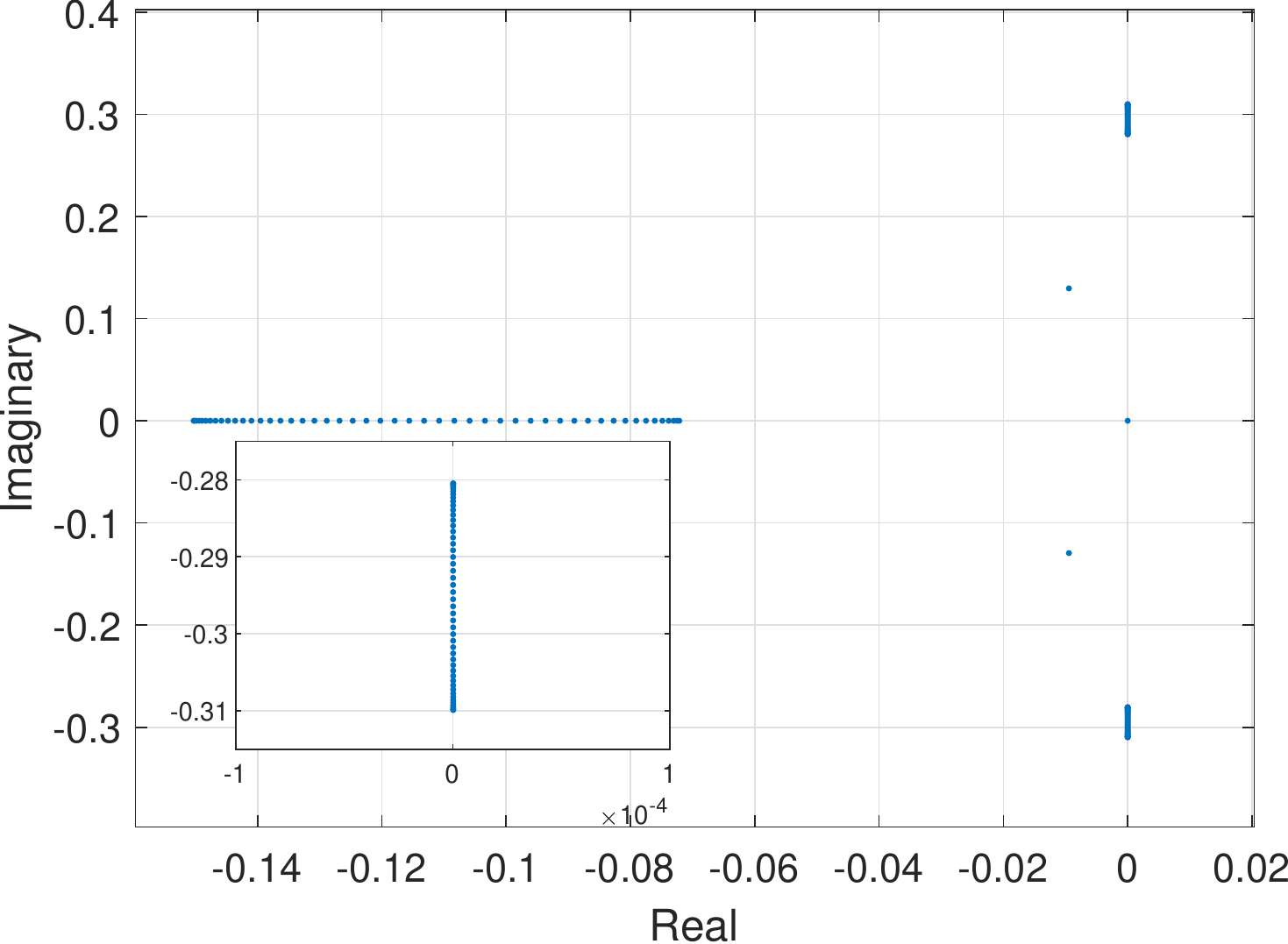}
\caption{Eigenvalues of the linearisation about a fixed point of~\eqref{eq:dadt}-\eqref{eq:dphidt}.
The inset shows a zoom of the main figure.
Parameters: $\mu=0.6,\nu=0.4,\alpha=\pi/2-0.08,
\Delta\omega=0.005$.}
\label{fig:eigskur}
\end{center}
\end{figure}

A typical set of eigenvalues of the linearisation about a fixed point 
of~\eqref{eq:dadt}-\eqref{eq:dphidt} is shown in Fig.~\ref{fig:eigskur}. Recall that we have
discretised $\omega$ with 50 points. We see 49 points on the negative real axis,
each corresponding to a perturbation localised at one or two
neighbouring $\phi$ values. There are also 49
complex conjugate pairs with zero real part, each corresponding to a perturbation
localised at one or two neighbouring 
$a$ values. These 147 eigenvalues are associated with discretising the continuous parameter $\omega$,
and are presumably
discretisations of the continuous spectrum associated with fixed points
of~\eqref{eq:dadt}-\eqref{eq:dphidt}.
There is also a complex conjugate pair with negative real 
part which, upon varying the relative sizes of $\mu$ and $\nu$, could cross the imaginary
axis resulting in a Hopf bifurcation~\cite{abrmir08}, and the single zero eigenvalue
corresponding to the invariance of the system under a global phase shift. 
As $\Delta\omega\to 0$, the 49 negative real eigenvalues collapse to a single negative real
value with multiplicity 49, and the 49 complex conjugate pairs on the imaginary axis collapse
to a single complex conjugate pair on the imaginary axis, again with multiplicity 49. 

Thus the term ``stable'' when referring to solutions in Fig.~\ref{fig:vardelomkur} actually
means ``neutrally stable'' or ``not unstable''. 
Indeed, when numerically integrating~\eqref{eq:dadt}-\eqref{eq:dphidt}
the system may not approach a fixed point even in a rotating coordinate frame, but the fixed
point can still be found using Newton's method. A similar phenomenon was observed in~\cite{pikros08},
and this is reflective of the fact that when~\eqref{eq:dadt} is discretised in $\omega$, the
equation for each $\omega$ describes the dynamics of a network of identical oscillators, 
whose dynamics is more
fully described by the equations derived using the Watanabe/Strogatz ansatz~\cite{watstr93,watstr94}.
(The Ott/Antonsen ansatz is a special case of the Watanabe/Strogatz ansatz, corresponding to
an infinite number of identical oscillators, with a uniform distribution of certain 
constants~\cite{pikros08}.)

Note that chimeras have been seen in similar networks with uniformly distributed 
frequencies~\cite{zhabi16}, but in the models studied there the coupling strength within
and between populations was a function of the level of synchrony within the populations,
not a constant, as here.

\subsection{Heterogeneous between-population coupling strengths}
\label{sec:varynu}

\begin{figure}
\begin{center}
\includegraphics[width=13cm]{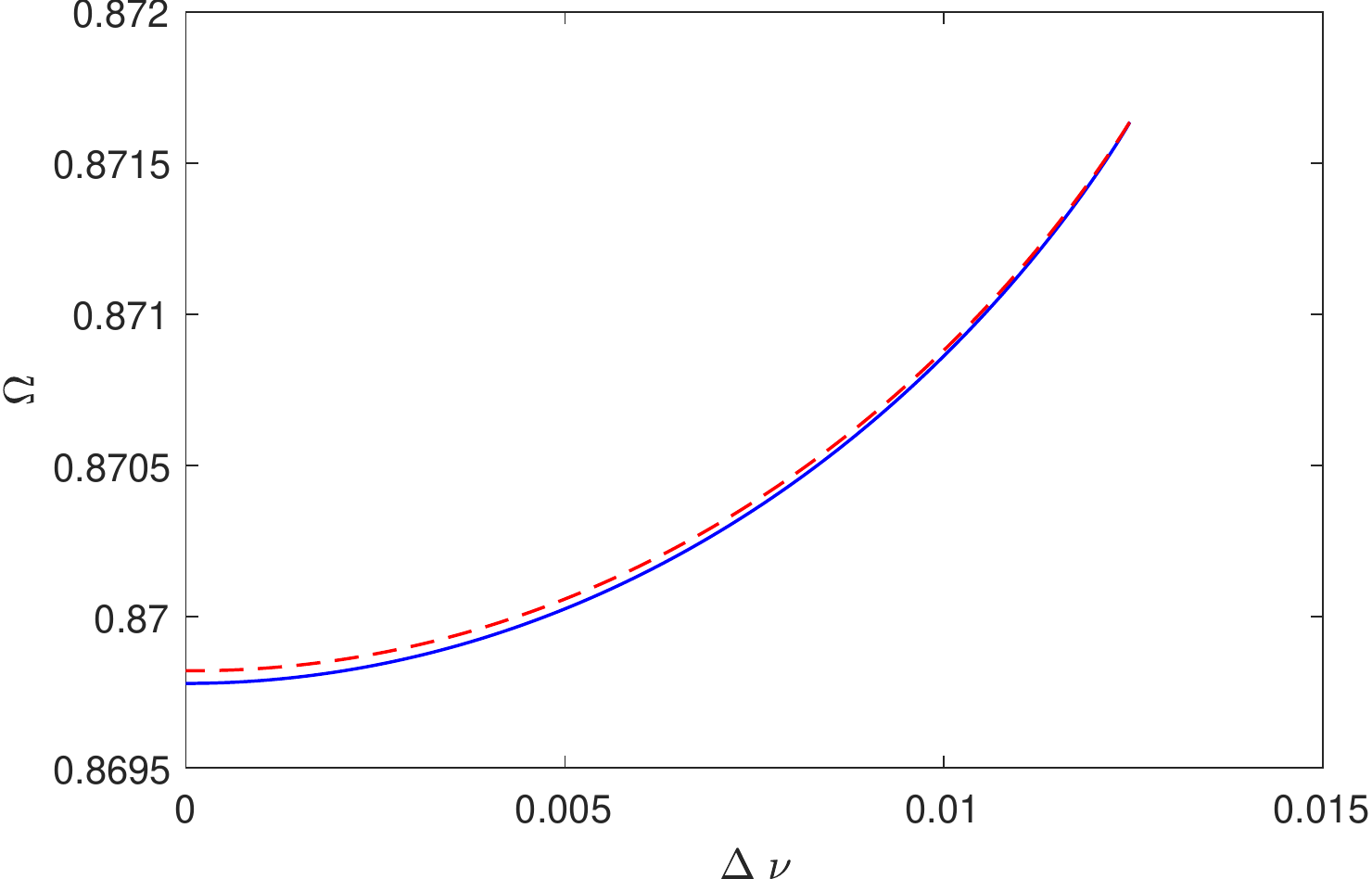}
\caption{$\Omega$ for 
fixed points of~\eqref{eq:dadt2}-\eqref{eq:dphidt2} describing chimera states. 
Solid: ``stable''. Dashed: unstable.
Parameters: $\mu=0.6,\nu_0=0.4,\alpha=\pi/2-0.08$.}
\label{fig:vardelnukur}
\end{center}
\end{figure}

Now, as another example, consider having heterogeneous 
$\nu$ values, i.e.~replace~\eqref{eq:dth1} with
\be
   \frac{d\theta_j}{dt}  = \omega+\frac{\mu}{N}\sum_{k=1}^N\sin{(\theta_k-\theta_j-\alpha)}+\frac{\nu_j}{N}\sum_{k=1}^N\sin{(\theta_{N+k}-\theta_j-\alpha)} \label{eq:dth1A}
\ee
and replace~\eqref{eq:dth2} in an equivalent way. Take the $\nu_j$ from $p(\nu)$,
a uniform distribution on $[\nu_0-\Delta\nu,\nu_0+\Delta\nu]$. In a similar
way to above we derive the evolution equations
\begin{align}
    \frac{\partial a(\nu,t)}{\partial t} & =-i\Omega a+\frac{1}{2}\left(e^{-i\alpha}Z-e^{i\alpha}\bar{Z}a^2\right) \label{eq:dadt2} \\
   \frac{\partial \phi(\nu,t)}{\partial t} & =-\Omega+(\mu\hat{S}+\nu S)\cos{(\phi+\alpha)}-(\mu\hat{C}+\nu C)\sin{(\phi+\alpha)} \label{eq:dphidt2}
\end{align}
where now
\be
   \hat{S}=\int_B \sin{[\phi(\nu,t)]}p(\nu)\ d\nu; \qquad \hat{C}=\int_B \cos{[\phi(\nu,t)]}p(\nu)\ d\nu
\ee
and
\be
    C+iS=\int_B a(\nu,t)p(\nu)\ d\nu,
\ee
$B$ is the interval $[\nu_0-\Delta\nu,\nu_0+\Delta\nu]$,
and without loss of generality
 we have set $\omega=0$. Note that $Z$, defined through~\eqref{eq:Z}, is not longer a scalar,
but a function of the continuous parameter $\nu$. Following fixed points 
of~\eqref{eq:dadt2}-\eqref{eq:dphidt2} as $\Delta\nu$ is increased
we obtain Fig.~\ref{fig:vardelnukur}, which
is very similar to Fig.~\ref{fig:vardelomkur}. The eigenvalues of
the linearisation about a ``stable'' state in Fig.~\ref{fig:vardelnukur}
are similar to those shown in Fig.~\ref{fig:eigskur}, for similar reasons as discussed above. 

We could perform similar analyses 
for the other two
parameters, $\alpha$ and $\mu$, but now move on to oscillators described by two variables.

%\cite{clupol16} for stability.

%van der Pol~\cite{omezak15}

%SNIC~\cite{vulhiz14}

%With mass~\cite{boukan14,olm15,belbri16}

%Two rings of SL~\cite{precha16}

%Two populations of SL~\cite{precha17}

%FHN bipartite~\cite{wuche18}

%LIF~\cite{olmpol11}

%QIF~\cite{ratpyr17}

%\cite{marbic16,monkur04,panabr16,pikros08}

\section{Stuart-Landau oscillators}
\label{sec:SL}
We now consider the chimera state found in~\cite{lai10} in a network of two populations of
Stuart-Landau oscillators, each oscillator being described by a complex variable. 

\subsection{Heterogeneous frequencies}
The
equations governing the dynamics are
\begin{align}
   \frac{dX_j}{dt} & =i\omega_j X_j+\epsilon^{-1}\{1-(1+\delta\epsilon i)|X_j|^2\}X_j+e^{-i\alpha}\left(\frac{\mu}{N}\sum_{k=1}^N X_k+\frac{\nu}{N}\sum_{k=1}^N X_{N+k}\right)
\label{eq:dXdt1}
\end{align}
for $j=1,\dots N$ and
\begin{align}
   \frac{dX_j}{dt} & =i\omega_j X_j+\epsilon^{-1}\{1-(1+\delta\epsilon i)|X_j|^2\}X_j +e^{-i\alpha}\left(\frac{\mu}{N}\sum_{k=1}^N X_{N+k}+\frac{\nu}{N}\sum_{k=1}^N X_{k}\right)
\label{eq:dXdt2}
\end{align}
for $j=N+1,\dots 2N$, where each $X_j\in\mathbb{C}$ and $\epsilon,\delta,\alpha,\mu$ and
$\nu$ are all real parameters. As before, $\mu$ is the strength of coupling within a 
population and $\nu$ is the strength between populations.
The $\omega_j$ are randomly chosen from the uniform distribution on
$[-\Delta\omega,\Delta\omega]$.

\begin{figure}
\begin{center}
\includegraphics[width=13cm]{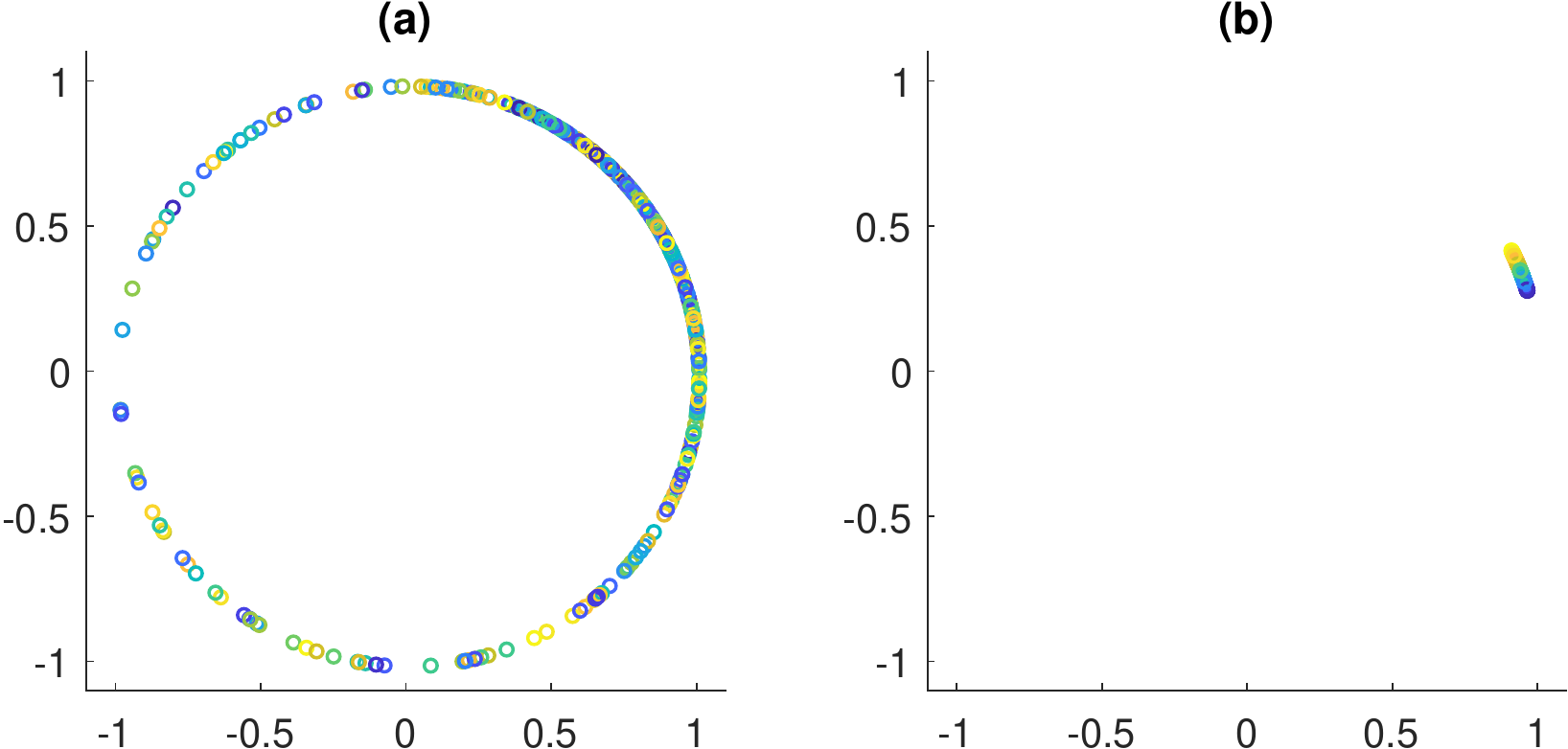}
\caption{A snapshot of a chimera state solution of~\eqref{eq:dXdt1}-\eqref{eq:dXdt2},
showing the $X_j$ in the complex plane. The points are coloured by their $\omega_j$ value.
(a): population 1; (b): population 2. 
Parameters: 
$N=500,\epsilon=0.05,\delta=-0.01,\mu=0.6,\nu=0.4,\alpha=\pi/2-0.08,\Delta\omega=0.01$.}
\label{fig:snapSL}
\end{center}
\end{figure}

An example of a stable chimera for~\eqref{eq:dXdt1}-\eqref{eq:dXdt2} 
is shown in Fig.~\ref{fig:snapSL}, with oscillators coloured by their $\omega_j$ value. 
(A similar figure appears in~\cite{lai19}.)
We see that population 2 is synchronised
and the oscillators lie on an open curve, with their position on the curve determined by their 
heterogeneous parameter $\omega_j$. Population 1 is incoherent, and there seems to be
no correlation between an oscillator's position and its $\omega_j$ value.
The oscillators in population 1 seem to all lie on a single closed curve, but we will see below
that oscillators with different values of $\omega_j$ actually lie on slightly different curves,
and move along these curves with slightly different average frequencies.

To analyse a chimera state let $X_{N+j}=Y_j$ for $j\in\{1,\dots N\}$ where the $Y_j$ rotate 
around the origin at the
same speed, i.e.~population two is
synchronised. Letting
\be
   \widehat{X}=\frac{1}{N}\sum_{k=1}^N X_k \qquad \mbox{and} \qquad  \widehat{Y}=\frac{1}{N}\sum_{k=1}^N Y_k
\ee
we have
\be
   \frac{dY_j}{dt}  =i\omega_{N+j} Y_j+\epsilon^{-1}\{1-(1+\delta\epsilon i)|Y_j|^2\}Y_j
+e^{-i\alpha}\left(\mu \widehat{Y}+\nu\widehat{X}\right) \label{eq:dYdtC}
\ee
for $j=1,\dots N$
and each oscillator in population one satisfies
\begin{align}
   \frac{dX_j}{dt} & =i\omega_j X_j+\epsilon^{-1}\{1-(1+\delta\epsilon i)|X_j|^2\}X_j +e^{-i\alpha}\left(\mu \widehat{X}+\nu \widehat{Y}\right),
\end{align}
for $j=1,\dots N$. Converting to polar coordinates for population 1
by writing $X_j=r_je^{i\phi_j}$ we have
\begin{align}
   \frac{dr_j}{dt} & = \epsilon^{-1}(1-r_j^2)r_j+\mbox{Re}\left[e^{-i(\alpha+\phi_j)}\left(\mu \widehat{X}+\nu \widehat{Y}\right)\right] \equiv F(r_j,\phi_j,\widehat{X},\widehat{Y}) \\
   \frac{d\phi_j}{dt} & = \omega_j-\delta r_j^2+\frac{1}{r_j}\mbox{Im}\left[e^{-i(\alpha+\phi_j)}\left(\mu \widehat{X}+\nu \widehat{Y}\right)\right]\equiv G(r_j,\phi_j,\widehat{X},\widehat{Y},\omega_j) \label{eq:G}
\end{align}
We now take the continuum limit of $N\to\infty$. Eqn~\eqref{eq:dYdtC} is replaced by
\be
    \frac{\partial Y}{\partial t}(\omega,t)  =i\omega Y+\epsilon^{-1}\{1-(1+\delta\epsilon i)|Y|^2\}Y+e^{-i\alpha}\left(\mu \widehat{Y}+\nu\widehat{X}\right) \label{eq:dYdt}
\ee
Generalising the theory from~\cite{clupol18,lai19}, we assume that oscillators
with a particular value of $\omega$ lie on a curve $\mathcal{C}(\omega)$ 
in the complex plane parametrised 
by the angle from the positive real axis, $\phi$. 
The distance from the origin to $\mathcal{C}(\omega)$ at angle $\phi$
is $R(\phi,t;\omega)$ and the density of oscillators at this point is $P(\phi,t;\omega)$.
The evolution of the functions $R$ and $P$ is given by
\begin{align}
   \frac{\p R}{\p t}(\phi,t;\omega) & = F(R,\phi,\widehat{X},\widehat{Y})-G(R,\phi,\widehat{X},\widehat{Y},\omega)\frac{\p R}{\p \phi} \label{eq:dRdt} \\
   \frac{\p P}{\p t}(\phi,t;\omega) & = -\frac{\p}{\p \phi}\left[P(\phi,t;\omega)G(R,\phi,\widehat{X},\widehat{Y},\omega)\right]+D\frac{\p^2}{\p \phi^2}P(\phi,t;\omega) \label{eq:dPdt}
\end{align}
where for numerical stability reasons we have added a small amount of diffusion, of strength $D$,
to~\eqref{eq:dPdt} (as did~\cite{clupol18}).
In the continuum limit we have
\be
   \widehat{X}=\int_B p(\omega)\int_0^{2\pi} P(\phi,t;\omega)R(\phi,t;\omega)e^{i\phi} d\phi \ d\omega \label{eq:Xb}
\ee
and
\be
   \widehat{Y}=\int_B Y(\omega,t)p(\omega)\ d\omega \label{eq:Yb}
\ee
where $B$ is the support of the unform density $p(\omega)$.
 The equations~\eqref{eq:dYdt}-\eqref{eq:Yb}
form a set of PDEs, coupled through integrals.
%We define $\beta=\pi/2-\alpha$.
Note that~\eqref{eq:dXdt1}-\eqref{eq:dXdt2} are invariant under the global phase shift
$X_j\mapsto X_je^{i\gamma}$ for any constant $\gamma$ and thus we can move to a rotating
coordinate frame in which $Y(\omega,t)$ is constant. Moving to a coordinate frame rotating with speed $\Omega$
has the effect of replacing the $\omega_j$ in~\eqref{eq:G} by $\omega_j-\Omega$ 
and the $\omega$ in~\eqref{eq:dYdt}  by
$\omega-\Omega$.

\begin{figure}
\begin{center}
\includegraphics[width=13cm]{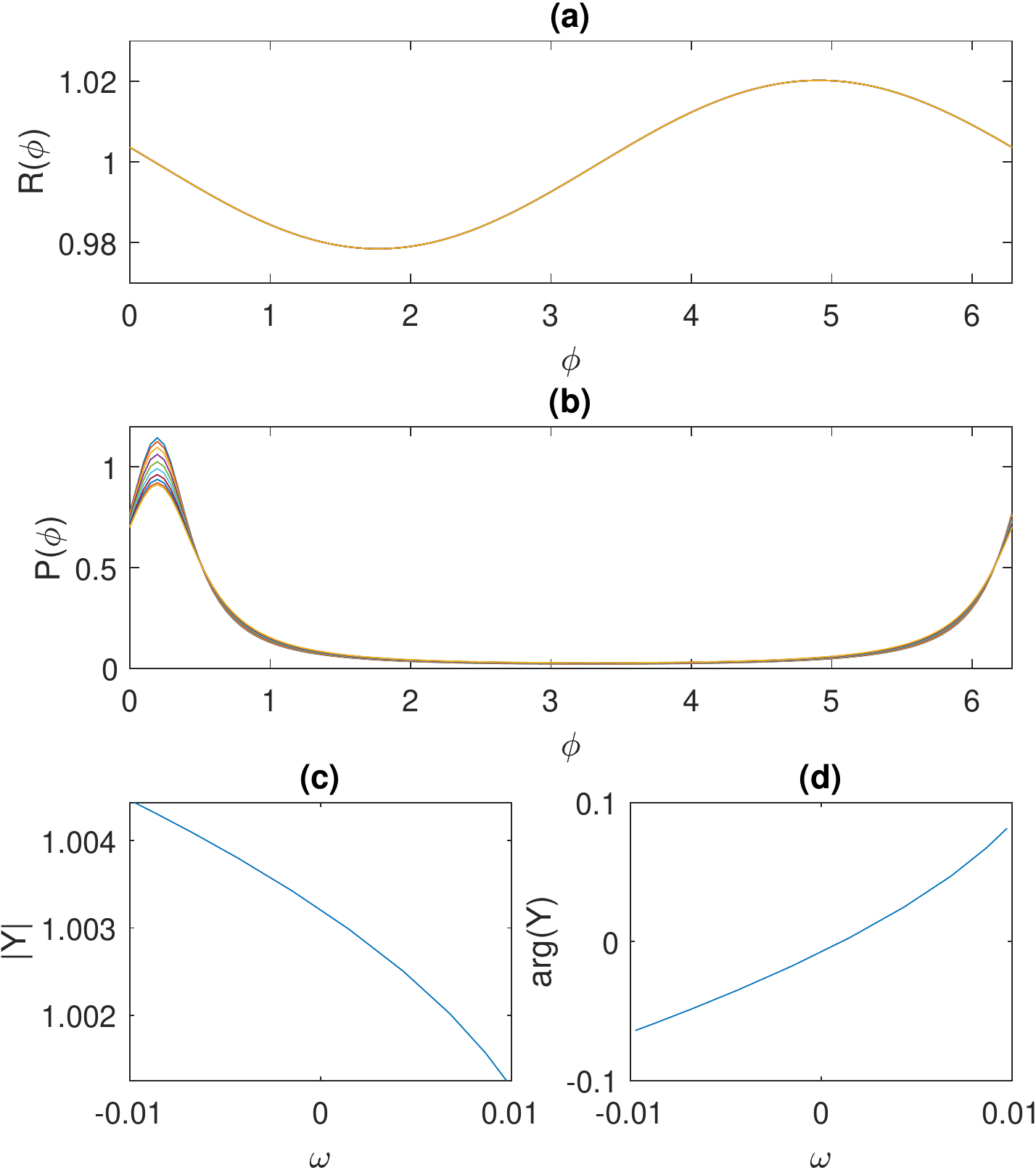}
\caption{Snapshot of a solution of~\eqref{eq:dRdt}-\eqref{eq:Yb} for which $\widehat{Y}$ is real.
This solution is stationary in a coordinate frame rotating at $\Omega=0.86678$.
(a): $R(\phi)$ for the 10 different values of $\omega_i$, (b): $P(\phi)$ for the 
10 different values of $\omega_i$. (c) and (d): modulus and argument of $Y$, respectively,
as functions of 
$\omega$.
Parameters: $\epsilon=0.05,\delta=-0.01,\mu=0.6,\nu=0.4,\alpha=\pi/2-0.08,D=10^{-8},
\Delta\omega=0.01$.}
\label{fig:ss}
\end{center}
\end{figure}

\subsubsection{Results}

We numerically integrate~\eqref{eq:dYdt}-\eqref{eq:Yb} in time to find a stable solution.
An example is shown in Fig.~\ref{fig:ss}. 
We see that $R(\phi)$  depends very weakly on $\omega$
(the distance between curves in panel (a) is $\sim 10^{-5}$) whereas $P(\phi)$ depends more
strongly on $\omega$.
We discretised $\phi$ using
$128$ equally-spaced points and implemented derivatives with respect to $\phi$
spectrally~\cite{tre00}. 
$p(\omega)$ is uniform on $[-\Delta\omega,\Delta\omega]$ and we
implement the integrals over $\omega$ in~\eqref{eq:Xb}-\eqref{eq:Yb}
using Gauss-Legendre quadrature with 10 points, 
so the discrete values of $\omega$ used are $\omega_i=\Delta\omega x_i$ for $i=1,2,\dots 10$
where the $x_i$ are the roots of $P_{10}(x)$, the  Legendre polynomial of order 10.
%The integrals over $\omega$ are thus approximated by weighted sums.
For each $\omega_i$ 
we enforce conservation of probability by setting $P$ at one angular
grid point equal to $1/(\Delta\phi)$ minus the sum of the values at all other grid points,
where $\Delta\phi=2\pi/128$, the $\phi$ grid spacing~\cite{erm06}.

\begin{figure}
\begin{center}
\includegraphics[width=13cm]{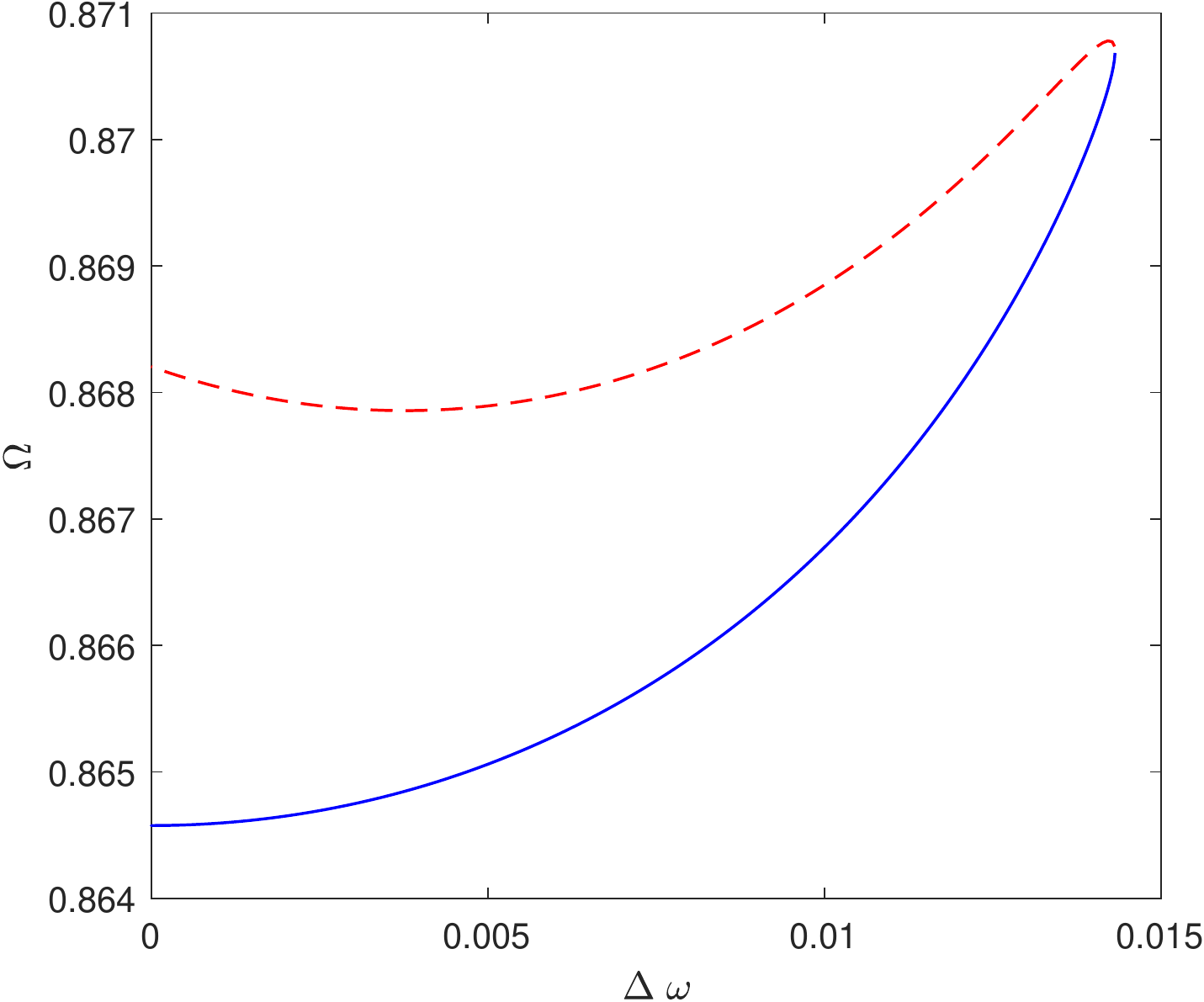}
\caption{$\Omega$ as a function of $\Delta\omega$ for fixed points (in a rotating
coordinate frame) of~\eqref{eq:dYdt}-\eqref{eq:Yb}
describing a chimera.
Solid: stable; dashed: unstable.
Parameters: $\epsilon=0.05,\delta=-0.01,\mu=0.6,\nu=0.4,\alpha=\pi/2-0.08,D=10^{-8}$.}
\label{fig:varydelom}
\end{center}
\end{figure}

Moving to a rotating coordinate frame and following a steady state 
of~\eqref{eq:dYdt}-\eqref{eq:Yb} in that frame as $\Delta\omega$ 
is increased we obtain Fig.~\ref{fig:varydelom}. As in Sec.~\ref{sec:varyomkur}
we see that the stable solution is destroyed in a saddle-node bifurcation.
The eigenvalues of the linearisation about a stable state are similar to those
in Fig.~3 of~\cite{lai19} and Fig.~5 in~\cite{clupol18}, i.e.~they form two clusters (not shown).
Those in the cluster 
with large negative real part are associated with perturbations in the $R$ component
of the dynamics while those in the other cluster 
which are almost marginally stable are associated with
perturbations in the $P$ component. 
%At the saddle-node bifurcation the $Y$ component of the eigenvector is localised at the
%most extreme value of $\omega_i$, indicating that in a simulation of a finite network, it
%would be the oscillator with the largest value of $\omega_j$ that would ``unlock'' first.

\begin{figure}
\begin{center}
\includegraphics[width=13cm]{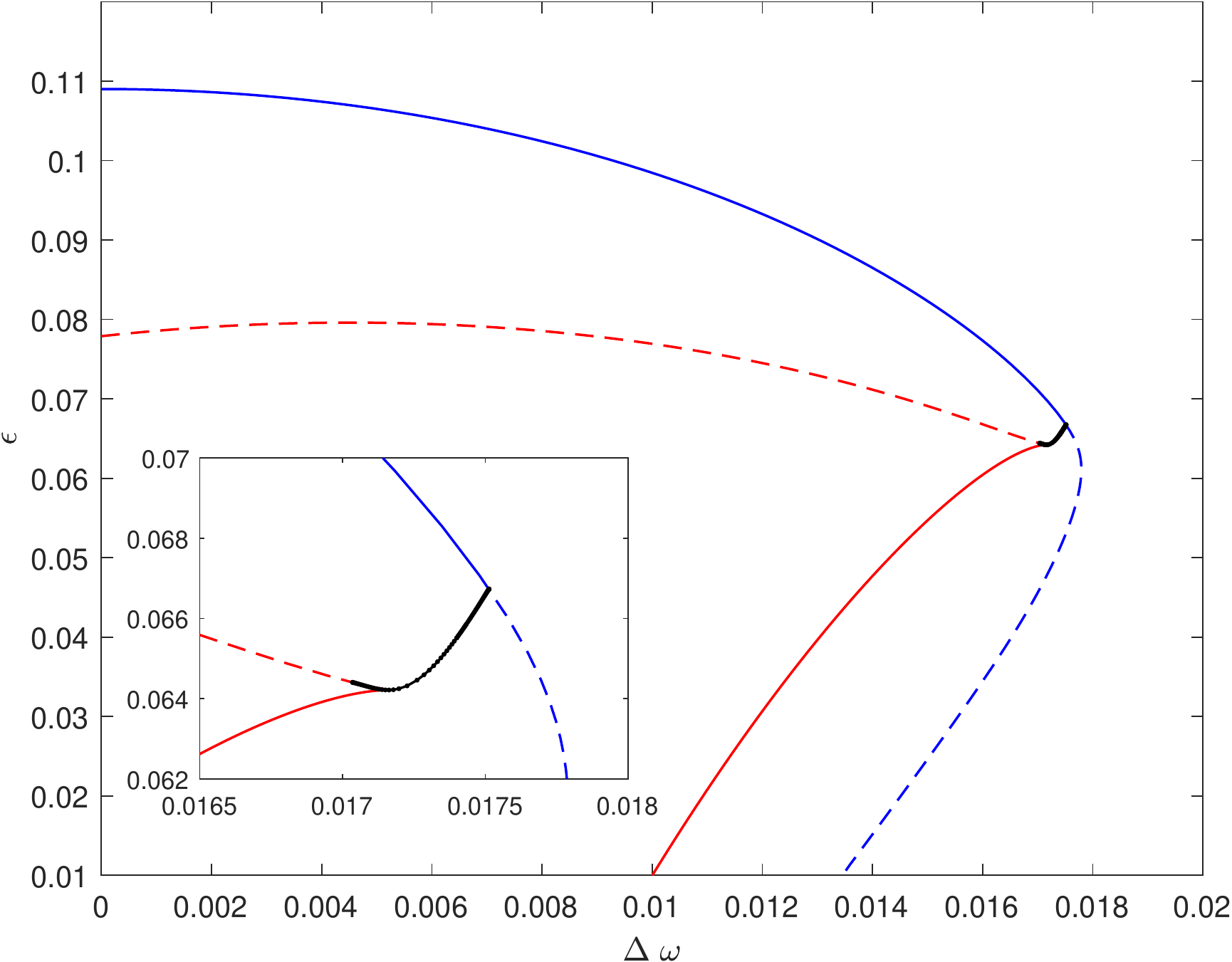}
\caption{Saddle-node bifurcation curves (red and blue) and Hopf bifurcation curve (black)
for fixed points of chimera solutions of~\eqref{eq:dYdt}-\eqref{eq:Yb}.
The inset shows a zoom of the Hopf bifurcation curve. The chimera is stable to the left of and below
the solid curves, and above the Hopf bifurcation curve.
Parameters: $\delta=-0.01,\mu=0.6,\nu=0.4,\alpha=\pi/2-0.08,D=10^{-8}$.}
\label{fig:epdelom}
\end{center}
\end{figure}

We know from~\cite{lai19} that if $\Delta\omega=0$, increasing $\epsilon$ will also destroy the
chimera in a saddle-node bifurcation. Following this bifurcation and the one in 
Fig.~\ref{fig:varydelom} we obtain Fig.~\ref{fig:epdelom}, where the former bifurcation is
shown in blue and the latter in red. Interestingly, they are not part of the same curve.
Both curves of saddle-node bifurcations have Takens-Bogdanov points on them, at which 
point the linearisation of the dynamics around the fixed point has a double zero
eigenvalue~\cite{guchol83,kuz04}. 
At both of these points a curve
of Hopf bifurcations is created, as seen in Fig.~\ref{fig:epdelom}. We expect other global
bifurcations in the vicinity of these points, but finding them numerically is difficult.

\subsection{Heterogeneous between-population coupling strengths}

\begin{figure}
\begin{center}
\includegraphics[width=13cm]{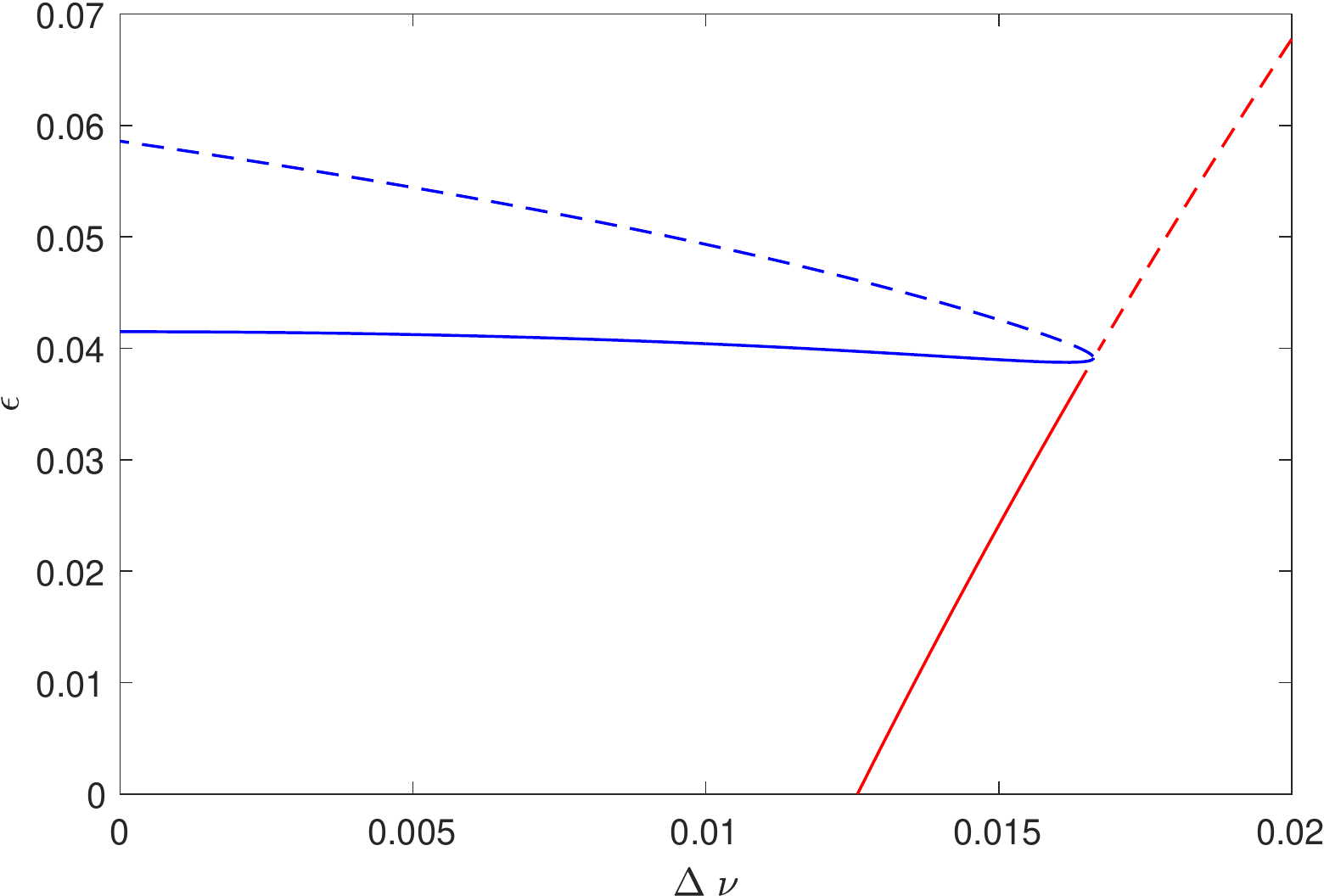}
\caption{Saddle-node bifurcation (red) and Hopf bifurcation (blue)
for fixed points of~\eqref{eq:dYdt}-\eqref{eq:Yb}. There is a stable stationary
chimera in the region bounded by the axes and the solid curves. There is a stable periodic
chimera above the lower (solid blue) Hopf bifurcation curve.
Parameters: $\delta=-0.01,\mu=0.625,\nu_0=0.375,\alpha=\pi-0.08,D=10^{-8},\omega=0$.}
\label{fig:vardelnuep}
\end{center}
\end{figure}

Now consider heterogeneity in $\nu$. We replace~\eqref{eq:dXdt1} by
\begin{align}
   \frac{dX_j}{dt} & =i\omega X_j+\epsilon^{-1}\{1-(1+\delta\epsilon i)|X_j|^2\}X_j+e^{-i\alpha}\left(\frac{\mu}{N}\sum_{k=1}^N X_k+\frac{\nu_j}{N}\sum_{k=1}^N X_{N+k}\right)
\end{align}
for $j=1,2,\dots N$ and similarly for population 2,
and choose the $\nu_j$ from a uniform distribution on 
$[\nu_0-\Delta\nu,\nu_0+\Delta\nu]$. Choosing $\mu=0.625$ and $\nu_0=0.375$ (and 
$\delta=-0.01,\alpha=\pi/2-0.08$) we know from~\cite{lai19} that if $\Delta\nu=0$, 
increasing $\epsilon$ results in
the stable chimera undergoing a supercritical Hopf bifurcation. Increasing $\Delta\nu$
for $\epsilon=0.03$ we find a saddle-node bifurcation at $\Delta\nu\approx 0.156$
and following these two bifurcations we obtain Fig.~\ref{fig:vardelnuep}.

The curves meet in a saddle-node/Hopf bifurcation, where the linearisation about the
fixed point has both a zero eigenvalue and a complex conjugate pair of purely imaginary
eigenvalues~\cite{guchol83,kuz04}. Stationary chimeras of the form we are
considering (where $Y$ is constant in a uniformly-rotating coordinate frame) exist only
to the left of the saddle-node bifurcation curve. Stable stationary states of this form only
exist in the region bounded by the axes and the solid curves in Fig.~\ref{fig:vardelnuep}.
They become unstable through a supercritical Hopf bifurcation as $\epsilon$ is increased,
leading to stable periodic chimeras. 

To better understand Fig.~\ref{fig:vardelnuep},
for $\epsilon=0.05$, as $\Delta\nu$ is decreased and the dashed saddle-node curve 
is crossed, a pair
of fixed points is created, one with three unstable eigenvalues and one with two.
As $\Delta\nu$ is further decreased the fixed point with three unstable eigenvalues
undergoes a Hopf bifurcation, gaining two stable directions. So between the solid and
dashed Hopf bifurcation curves one fixed point has one unstable eigenvalue and the
other has two. As $\epsilon$ is then decreased the fixed point with two unstable eigenvalues
undergoes a Hopf bifurcation, becoming stable. If $\Delta\nu$ is then increased, this stable
fixed point is destroyed in a saddle-node bifurcation with the fixed point having one
unstable eigenvalue. We expect there to be other curves of global bifurcations
in a neighbourhood of the saddle-node/Hopf bifurcation, but finding them is numerically
difficult.

\section{Kuramoto with inertia}
\label{sec:pend}

We now consider a network formed from two populations of $N$ Kuramoto oscillators with inertia,
where we have heterogeneity in frequencies. 
The system is described by
\begin{align}
   m\frac{d^2\theta_i^{(1)}}{dt^2}+\frac{d\theta_i^{(1)}}{dt} & =\omega_i^{(1)}+\frac{\mu}{N}\sum_{j=1}^N\sin{\left(\theta_j^{(1)}-\theta_i^{(1)}-\alpha\right)}+\frac{\nu}{N}\sum_{j=1}^N\sin{\left(\theta_j^{(2)}-\theta_i^{(1)}-\alpha\right)} \label{eq:pendA} \\
   m\frac{d^2\theta_i^{(2)}}{dt^2}+\frac{d\theta_i^{(2)}}{dt} & =\omega_i^{(2)}+\frac{\mu}{N}\sum_{j=1}^N\sin{\left(\theta_j^{(2)}-\theta_i^{(2)}-\alpha\right)}+\frac{\nu}{N}\sum_{j=1}^N\sin{\left(\theta_j^{(1)}-\theta_i^{(2)}-\alpha\right)} \label{eq:pendB}
\end{align}
where $m$ is ``mass'', $\mu,\nu$ and $\alpha$ are parameters,
and the superscript labels the population. When $m=0$ and the $\omega_i$ are all equal
this reverts to a previously studied
case~\cite{abrmir08,pikros08}. With $m=0$ and $\omega_i$ chosen from a uniform distribution
it reverts to that studied in Sec.~\ref{sec:varomkur}, while for $m\neq 0$ 
and the $\omega_j$ being chosen
from a Lorentzian it is that studied in~\cite{boukan14}. These authors
found apparently stable chimeras for finite networks of oscillators.
With $m\neq 0$ and identical
$\omega_j$ it is the same as studied in~\cite{olm15}.

In~\cite{lai19} it was found that with $m \neq 0$ and
identical $\omega_i$ the apparently stable chimera solution found by 
simulating~\eqref{eq:pendA}-\eqref{eq:pendB}
was actually (weakly) unstable when the continuum equations were studied, with
the real part of the rightmost eigenvalues determining its stability increasing with $m$. 
So it is of interest to investigate the effects of heterogeneity in the $\omega_j$
on such a state: does it stabilise this state?

\begin{figure}
\begin{center}
\includegraphics[width=13cm]{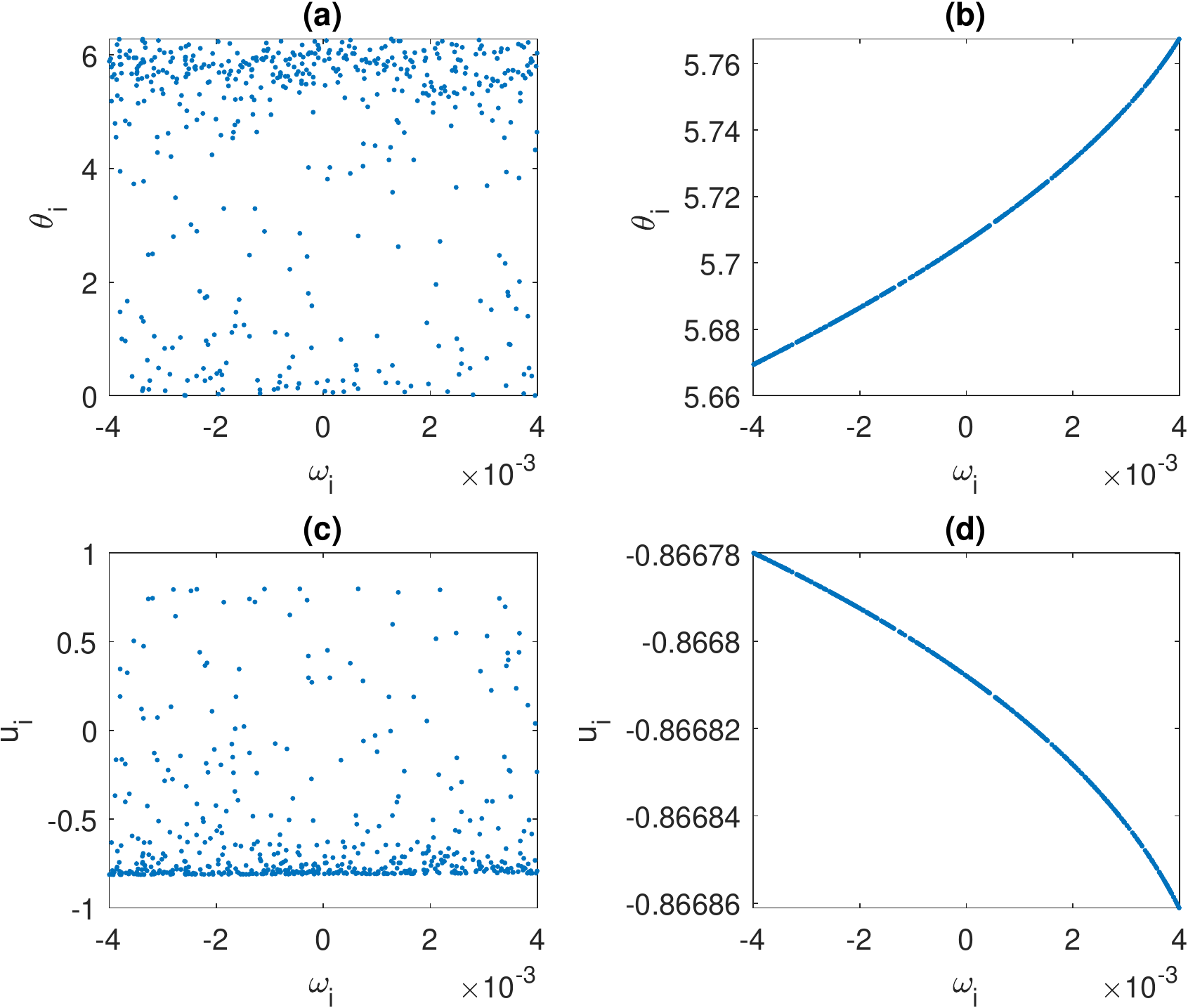}
\caption{A snapshot of a chimera state for~\eqref{eq:pendAA}-\eqref{eq:pendDD}. 
Population 1 is on the left and population 2 is on the right. Note the different vertical scales.
Parameters: $N=500,m=0.1,\mu=0.6,\nu=0.4,\alpha=\pi/2-0.05,\Delta\omega=0.004$.}
\label{fig:pendsnap}
\end{center}
\end{figure}

%We note that this system is invariant under a uniform shift
%of all of the phases.
We rewrite the equations as
\begin{align}
   \frac{d\theta_i^{(1)}}{dt} & = u_i^{(1)} \label{eq:pendAA} \\
   \frac{du_i^{(1)}}{dt} & = \left[\omega_i^{(1)}-u_i^{(1)}+\frac{\mu}{N}\sum_{j=1}^N\sin{\left(\theta_j^{(1)}-\theta_i^{(1)}-\alpha\right)} \right. \left. +\frac{\nu}{N}\sum_{j=1}^N\sin{\left(\theta_j^{(2)}-\theta_i^{(1)}-\alpha\right)}\right]/m \\
    \frac{d\theta_i^{(2)}}{dt} & = u_i^{(2)} \\
   \frac{du_i^{(2)}}{dt} & = \left[\omega_i^{(2)}-u_i^{(2)}+\frac{\mu}{N}\sum_{j=1}^N\sin{\left(\theta_j^{(2)}-\theta_i^{(2)}-\alpha\right)} \right. \left. +\frac{\nu}{N}\sum_{j=1}^N\sin{\left(\theta_j^{(1)}-\theta_i^{(2)}-\alpha\right)}\right]/m \label{eq:pendDD}
\end{align}
A snapshot of a stable chimera state for~\eqref{eq:pendAA}-\eqref{eq:pendDD} is shown in 
Fig.~\ref{fig:pendsnap}, where for both populations the $\omega_i$ 
are taken from a uniform distribution
on $[-\Delta\omega,\Delta\omega]$. We see that population 1 is incoherent while population 2 is
synchronised with both $\theta_i$ and $u_i$ being smooth functions of $\omega_i$.

To analyse this state let us assume that
population two is synchronised, with $\theta_k^{(2)}=\Theta_k$ and $u_k^{(2)}=U_k$ 
for $k=1,2\dots N$. We drop the superscripts for variables in population 1.
%Due to the invariance under global phase shift we can move to a coordinate frame rotating at 
%constant speed $\Omega$ 
Oscillators in population 2 satisfy
\begin{align}
   \frac{d\Theta_k}{dt} & = U_k  \label{eq:dThetadt} \\
   \frac{dU_k}{dt} & = %\left[\omega_k^{(2)}-U_k+\mu\mbox{Im}\left\{e^{-i(\Theta_k+\alpha)}Y\right\}+\frac{\nu}{N}\sum_{j=1}^N\sin{(\theta_j-\Theta_k-\alpha)}\right]/m \nonumber \\
    \left[\omega_k^{(2)}-U_k+\mu\mbox{Im}\left\{e^{-i(\Theta_k+\alpha)}Y\right\}+\nu\mbox{Im}\left\{e^{-i(\Theta_k+\alpha)}X\right\}\right]/m \label{eq:dUdtA}
\end{align}
where
\be
   X\equiv\frac{1}{N}\sum_{j=1}^N e^{i\theta_j}\in\mathbb{C}, \qquad Y=\frac{1}{N}\sum_{j=1}^N e^{i\Theta_j}\in\mathbb{C}
\ee
 Oscillators in population 1 satisfy
\begin{align}
   \frac{d\theta_k}{dt} & = u_k \\
   \frac{du_k}{dt} & = %\left[\omega_k^{(1)}-u_k+\frac{\mu}{N}\sum_{j=1}^N\sin{(\theta_j-\theta_k-\alpha)}+\nu\mbox{Im}\left\{e^{-i(\theta_k+\alpha)}Y\right\}\right]/m \nonumber \\
     \left[\omega_k^{(1)}-u_k+\mu\mbox{Im}\left\{e^{-i(\theta_k+\alpha)}X\right\}+\nu\mbox{Im}\left\{e^{-i(\theta_k+\alpha)}Y\right\}\right]/m
\end{align}
for $k=1,\dots N$.
We put these equations in ``polar'' form by defining $r_k=2+u_k$ (adding 2 bounds the $r_k$
away from zero) and thus we have
\begin{align}
   \frac{dr_k}{dt} & = \left[\omega_k^{(1)}-(r_k-2)+\mu\mbox{Im}\left\{e^{-i(\theta_k+\alpha)}X\right\}  +\nu\mbox{Im}\left\{e^{-i(\theta_k+\alpha)}Y\right\}\right]/m \nonumber \\
   & \equiv F(r_k,\theta_k,X,Y,\omega_k^{(1)}) \\
   \frac{d\theta_k}{dt} & = r_k-2 \label{eq:dthetadt}
\end{align} 
For the chimera state of interest, $\Theta_k$ and $U_k$ are
stationary in a coordinate frame rotating at speed $\Omega$.
Moving to this coordinate frame has the effect of replacing~\eqref{eq:dThetadt} by 
\be
   \frac{d\Theta_k}{dt} = U_k+\Omega \label{eq:dThetadtA}
\ee
and~\eqref{eq:dthetadt} by
\be
   \frac{d\theta_k}{dt}= r_k-2+\Omega \equiv G(r_k,\theta_k,X,Y)
\ee
Taking the limit $N\rightarrow\infty$ we consider the dynamical system
\begin{align}
   \frac{\p R}{\p t}(\theta,t;\omega) & = F(R,\theta,X,Y,\omega)-G(R,\theta,X,Y)\frac{\p R}{\p \theta} \label{eq:dRpend} \\
   \frac{\p P}{\p t}(\theta,t;\omega) & = -\frac{\p}{\p \theta}\left[P(\theta,t;\omega)G(R,\theta,X,Y)\right]+D\frac{\p^2}{\p \theta^2}P(\theta,t;\omega) \label{eq:dPpend}
\end{align}
along with
\be
   \frac{\partial \Theta}{\partial t}(\omega,t) = U+\Omega \label{eq:dThetadtB}
\ee
and
\be
\frac{\partial U}{\partial t}(\omega,t) =\left[\omega-U+\mu\mbox{Im}\left\{e^{-i(\Theta+\alpha)}Y\right\}+\nu\mbox{Im}\left\{e^{-i(\Theta+\alpha)}X\right\}\right]/m 
\ee
where
\be
   X(t)=\int_B p(\omega) \int_0^{2\pi}P(\theta,t;\omega)R(\theta,t;\omega)e^{i\theta}d\theta\ d\omega \label{eq:X}
\ee
and
\be
   Y(t)=\int_B e^{i\Theta(\omega,t)}p(\omega) d\omega \label{eq:Ypend}
\ee
and $B$ is the interval $[-\Delta\omega,\Delta\omega]$.

\subsection{Results}

\begin{figure}
\begin{center}
\includegraphics[width=13cm]{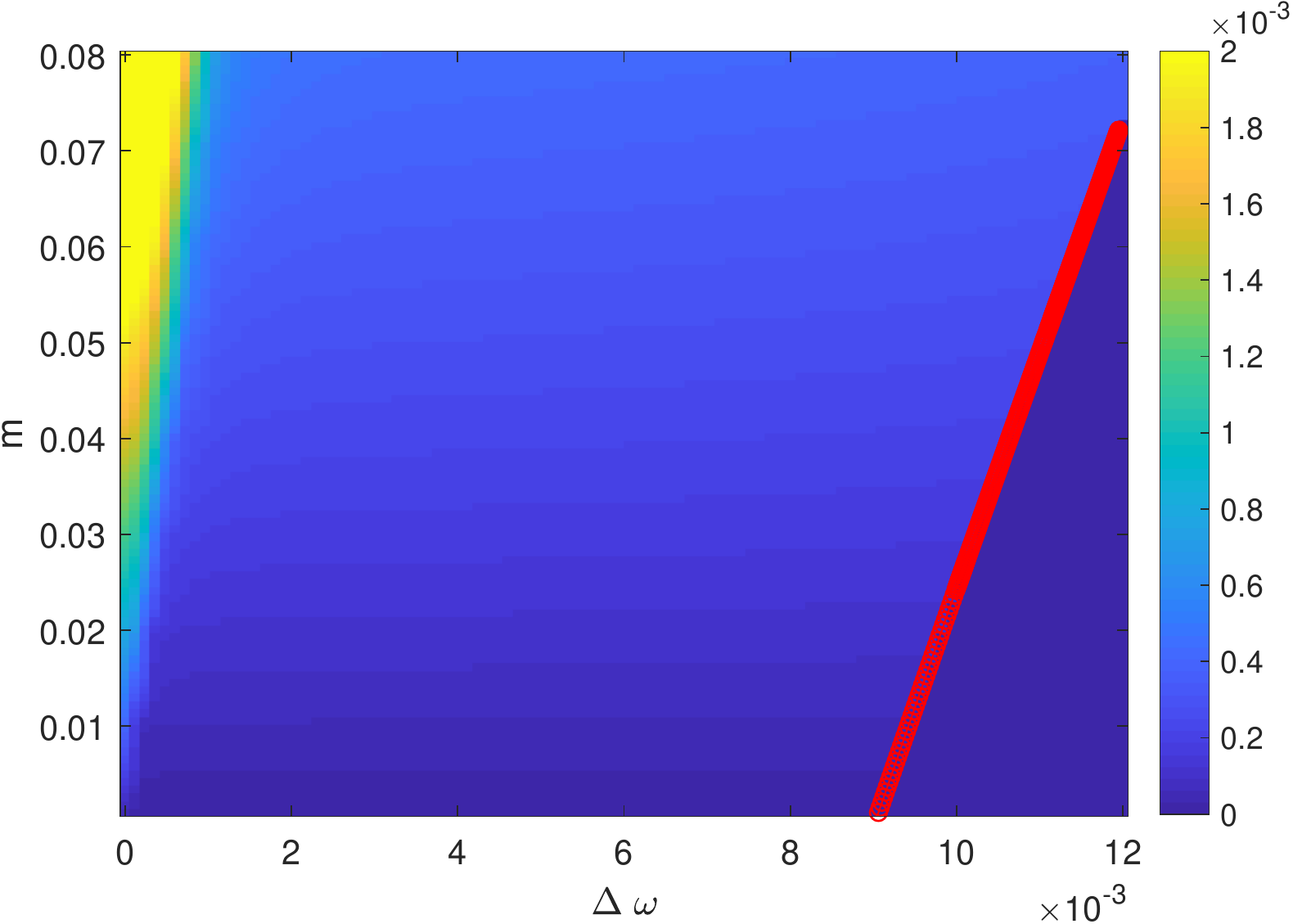}
\caption{The red circles show a saddle-node bifurcation, to the right of which solutions
in which all oscillators in population 2 for~\eqref{eq:dRpend}-\eqref{eq:Ypend} 
are locked do not exist. The colour is the real part
of the rightmost eigenvalues determining the stability of the chimera 
(truncated above $2\times 10^{-3}$), which is always positive.
Parameters: $\mu=0.6,\nu=0.4,\alpha=\pi/2-0.05,D=10^{-13}$. $\theta$ is discretised with 128
points and $\omega$ with 10.}
\label{fig:mdelompend}
\end{center}
\end{figure}

For small values of $m,\Delta\omega$ and $D$, stable chimera solutions 
of~\eqref{eq:dRpend}-\eqref{eq:Ypend} can be found, but as in~\cite{lai19}, decreasing
$D$ to $10^{-13}$ we find that they are actually weakly unstable. Following them as $\Delta\omega$
is increased we find that they are destroyed in a saddle-node bifurcation, as shown in
Fig.~\ref{fig:mdelompend}. To the left of this curve these solutions are always unstable,
although (as in~\cite{lai19}) they become less unstable as $m$ is decreased. Although
making the oscillators heterogeneous in this way does not fully stabilise the chimera,
it does make then less unstable, as can be seen by varying $\Delta\omega$ for a fixed value
of $m$ in Fig.~\ref{fig:mdelompend}. 
%This is consistent with the results of~\cite{boukan14}
%who found apparently stable chimeras for finite networks when the values of the $\omega_i$
%were chosen from a Lorentzian distribution.
In summary, as in~\cite{lai19}, for the parameters chosen, stable chimeras do not exist
for infinite networks described by~\eqref{eq:dRpend}-\eqref{eq:Ypend} even with $\omega$
values chosen from a uniform distribution.

\section{van der Pol oscillators}
\label{sec:vdp}
\begin{figure}
\begin{center}
\includegraphics[width=13cm]{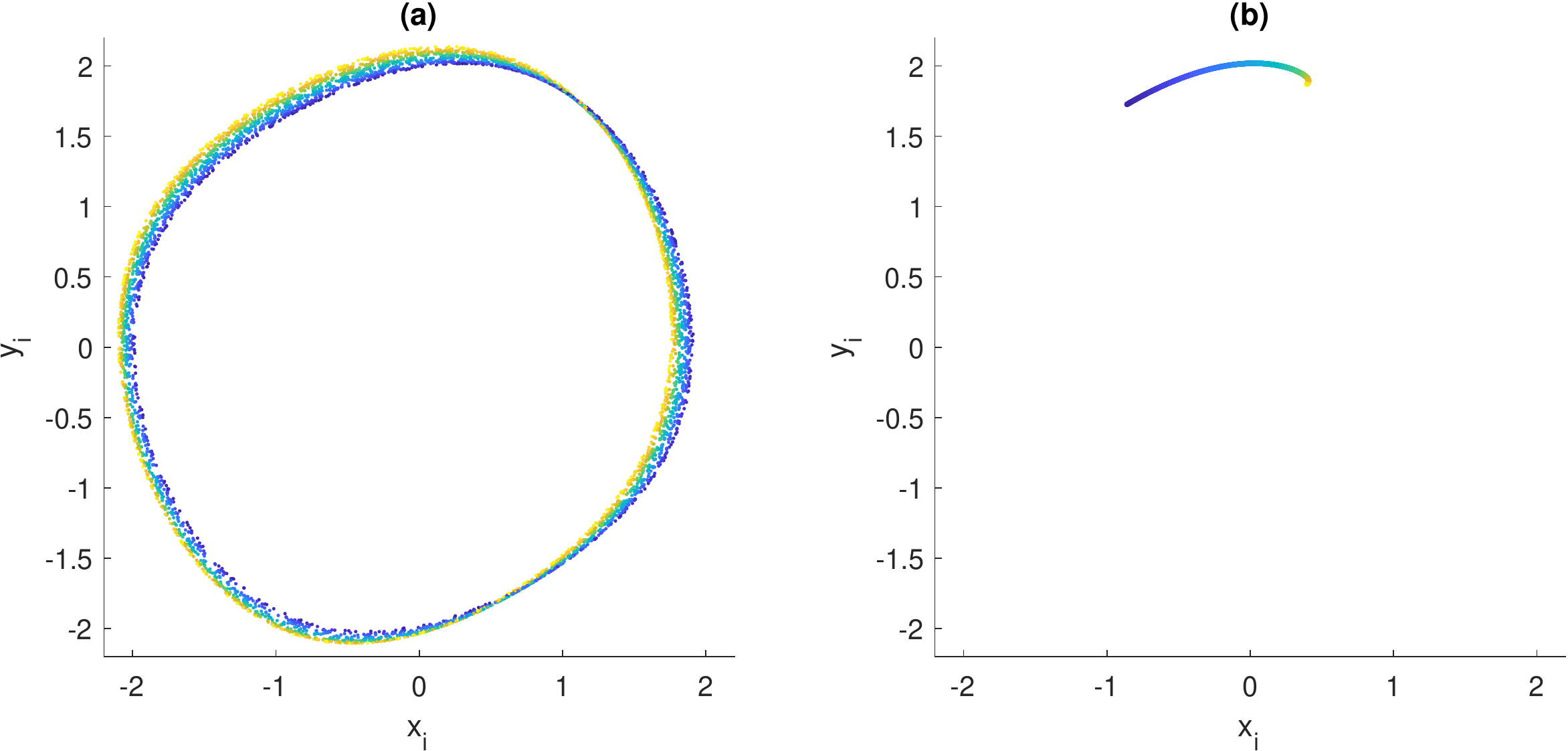}
\caption{Snapshot of a stable chimera solution of~\eqref{eq:dxdtvdp}-\eqref{eq:Mvdp}.
The oscillators are coloured by their $b_i$ value (blue low, yellow high).
(a): population 1; (b): population 2.
Parameters: $\mu=0.09,\nu=0.01,\epsilon=0.2,c=0.1,b_0=1,\Delta b=0.7,N=5000$.}
\label{fig:vdpsnap}
\end{center}
\end{figure}

We lastly consider two populations of van der Pol oscillators, governed by

%
%\subsection{van der Pol}
%
%
\begin{align}
   \frac{dx_i}{dt} & = y_i \label{eq:dxdtvdp} \\
   \frac{dy_i}{dt} & = \epsilon(1-x_i^2)y_i-x_i+\mu[b_i(\bar{X}_1-x_i)+c(\bar{Y}_1-y_i)]+\nu[b_i(\bar{X}_2-x_i)+c(\bar{Y}_2-y_i)] \label{eq:dydtvdp}
\end{align}
for $i=1,2\dots N$ and
\begin{align}
   \frac{dx_i}{dt} & = y_i \\
   \frac{dy_i}{dt} & = \epsilon(1-x_i^2)y_i-x_i+\mu[b_i(\bar{X}_2-x_i)+c(\bar{Y}_2-y_i)]+\nu[b_i(\bar{X}_1-x_i)+c(\bar{Y}_1-y_i)]
\end{align}
for $i=N+1,\dots 2N$, where
\be 
   \bar{X}_1=\frac{1}{N}\sum_{i=1}^N x_i; \qquad \bar{Y}_1=\frac{1}{N}\sum_{i=1}^N y_i;
\qquad \bar{X}_2=\frac{1}{N}\sum_{i=1}^N x_{N+i}; \qquad \bar{Y}_2=\frac{1}{N}\sum_{i=1}^N y_{N+i}
\label{eq:Mvdp}
\ee
As usual, $\mu$ is the within-population coupling strength and $\nu$ is the between-population
strength.  There is mean-field coupling involving both $x$ and $y$ variables.
Equations of this form were considered in~\cite{omezak15,uloome16} although with
nonlocal coupling on a ring of oscillators.  
We set $\epsilon=0.2,\mu=0.09,\nu=0.01,c=0.1$ and consider heterogeneity in the $b_i$,
choosing them from a uniform distribution on $[b_0-\Delta b,b_0+\Delta b]$. An example
of a stable chimera for $b_0=1,\Delta b=0.7$ is shown in Fig.~\ref{fig:vdpsnap}.
Population 1 is incoherent while population 2 is synchronised. For population 1 it is clear
that oscillators with different $b_i$ lie on different curves.

To analyse this state suppose 
population 2 is synchronised and $x_{N+i}=X_i$ and $y_{N+i}=Y_i$ for $i=1,2\dots N$.
These are the values shown in  Fig.~\ref{fig:vdpsnap}(b).
Then we have
\begin{align}
   \frac{dX_i}{dt} & = Y_i  \\
   \frac{dY_i}{dt} & = \epsilon(1-X_i^2)Y_i-X_i+\mu[b_i(\bar{X}_2-X_i)+c(\bar{Y}_2-Y_i)]+\nu[b_i(\bar{X}_1-X_i)+c(\bar{Y}_1-Y_i)]
\end{align}
For population 1, we see from Fig.~\ref{fig:vdpsnap}(a) that oscillators lie on curves
which completely contain the origin, so moving to polar coordinates by
writing $r_i^2=x_i^2+y_i^2$ and $\tan{\theta_i}=y_i/x_i$ so that
$x_i=r_i\cos{\theta_i}$ and $y_i=r_i\sin{\theta_i}$ we have
\begin{align}
	   \frac{dr_i}{dt} & =\frac{x_i\frac{dx_i}{dt}+y_i\frac{dy_i}{dt}}{r_i}\equiv F(r_i,\theta_i,\bar{X}_1,\bar{Y}_1,\bar{X}_2,\bar{Y}_2,b_i) \\
    \frac{d\theta_i}{dt} & = \frac{x_i\frac{dy_i}{dt}-y_i\frac{dx_i}{dt}}{r_i^2}\equiv G(r_i,\theta_i,\bar{X}_1,\bar{Y}_1,\bar{X}_2,\bar{Y}_2,b_i)
\end{align}
where $dx_i/dt$ and $dy_i/dt$ are given by~\eqref{eq:dxdtvdp}-\eqref{eq:dydtvdp}.
Taking the continuum limit 
we consider the dynamical system
\begin{align}
   \frac{\p R}{\p t}(\theta,t;b) & = F(R,\theta,\bar{X}_1,\bar{Y}_1,\bar{X}_2,\bar{Y}_2,b)-G(R,\theta,\bar{X}_1,\bar{Y}_1,\bar{X}_2,\bar{Y}_2,b)\frac{\p R}{\p \theta} \label{eq:dRvdp} \\
   \frac{\p P}{\p t}(\theta,t;b) & = -\frac{\p}{\p \theta}\left[P(\theta,t)G(R,\theta,\bar{X}_1,\bar{Y}_1,\bar{X}_2,\bar{Y}_2,b)\right]+D\frac{\p^2}{\p \theta^2}P(\theta,t;b) 
\end{align}
together with
\begin{align}
   \frac{\partial X(b,t)}{\partial t} & = Y \label{eq:dXdt} \\
   \frac{\partial Y(b,t)}{\partial t} & = \epsilon(1-X^2)Y-X+\mu[b(\bar{X}_2-X)+c(\bar{Y}_2-Y)]+\nu[b(\bar{X}_1-X)+c(\bar{Y}_1-Y)] \label{eq:dYdtD}
\end{align}
where
\be
   \bar{X}_1=\int_B\int_0^{2\pi}P(\theta,t;b)R(\theta,t;b)\cos{\theta}\ d\theta\ db; \qquad \bar{Y}_1=\int_B\int_0^{2\pi}P(\theta,t;b)R(\theta,t;b)\sin{\theta}\ d\theta \ db
\ee
and
\be
   \bar{X}_2=\int_B X(b,t)\ db; \qquad \bar{Y}_2=\int_B Y(b,t)\ db \label{eq:XYvdp}
\ee
where $B$ is the interval $[b_0-\Delta b,b_0+\Delta b]$.

\begin{figure}
\begin{center}
\includegraphics[width=13cm]{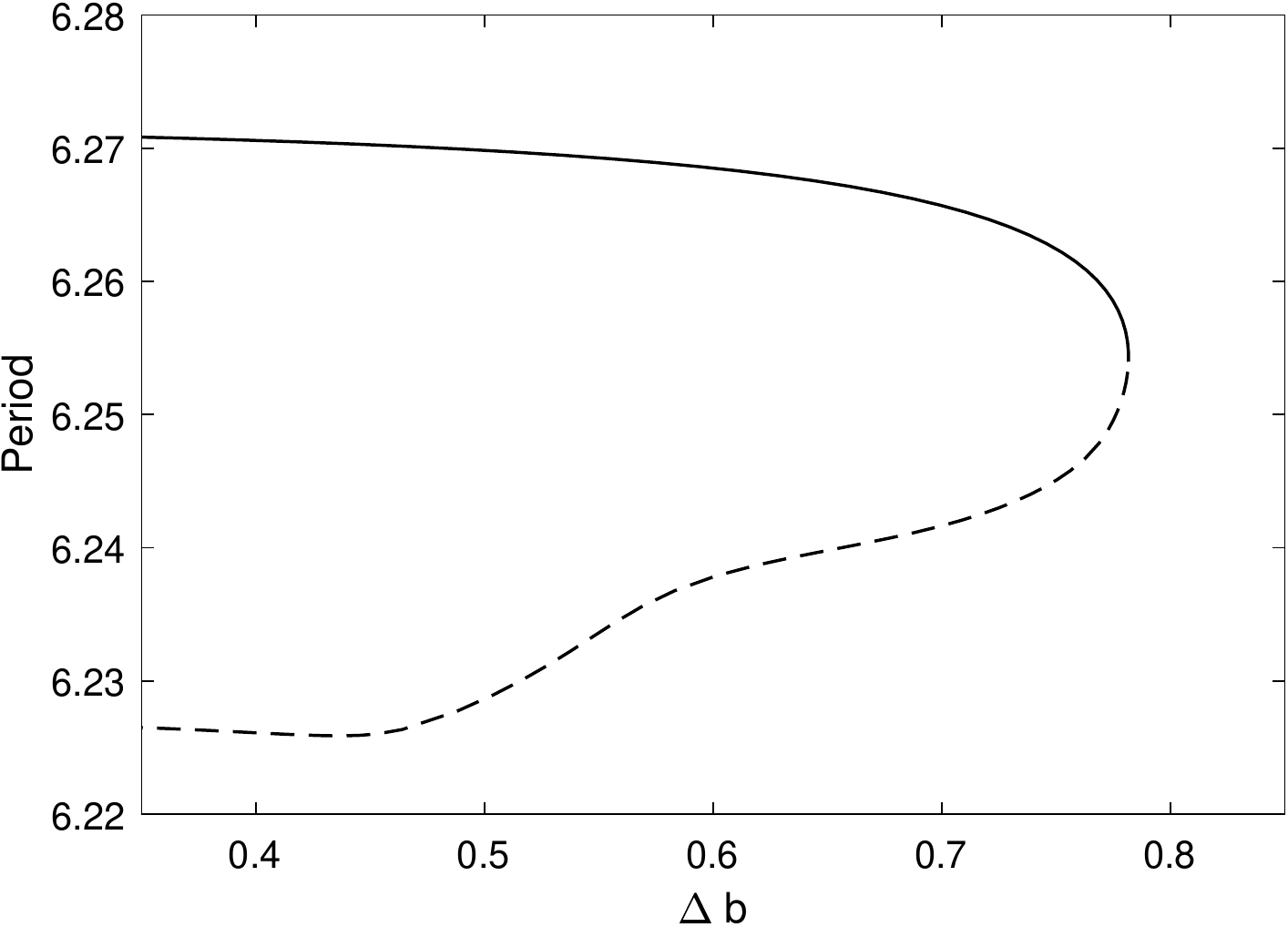}
\caption{Period of the periodic chimera solution of~\eqref{eq:dRvdp}-\eqref{eq:XYvdp}.
Solid: stable; dashed: unstable.
Parameters: $\mu=0.09,\nu=0.01,\epsilon=0.2,c=0.1,b_0=1,D=10^{-4}$.
$\theta$ is discretised in 128 points and $b$ in 10.}
\label{fig:vardbvdp}
\end{center}
\end{figure}

A significant difference between this system and those in the previous sections is that we 
cannot go a uniformly rotating coordinate frame in which $X$ and $Y$ are constant. Thus
the chimera of interest is a periodic solution of~\eqref{eq:dRvdp}-\eqref{eq:XYvdp}.
We numerically continue periodic solutions using pseudo-arclength continuation
and determine their stability in terms of the
magnitude of the Floquet multipliers of that solution. We obtain Fig.~\ref{fig:vardbvdp},
where we see the stable periodic solution is destroyed in a saddle-node bifurcation as $\Delta b$
is increased. Unlike the systems studied in Secs.~\ref{sec:kur} and~\ref{sec:SL}, the 
stable chimera
in this system is genuinely stable, without marginal or nearly marginal Floquet multipliers.

\section{Discussion}
\label{sec:disc}

We have considered chimeras in networks formed from two coupled populations of oscillators.
In each network one parameter has been chosen randomly from a uniform distribution. For narrow
enough distributions the chimeras which exist for the case of identical parameters persist,
and the synchronous oscillators remain synchronised, although no longer having identical states.
We generalised the theory in~\cite{lai19} to cover these states, at the price of increased
computational effort. In all cases we found that chimeras were destroyed in saddle-node
bifurcations as the width of the uniform distribution was increased. We now discuss
several generalisations of the approach taken here.

\begin{figure}
\begin{center}
\includegraphics[width=15cm]{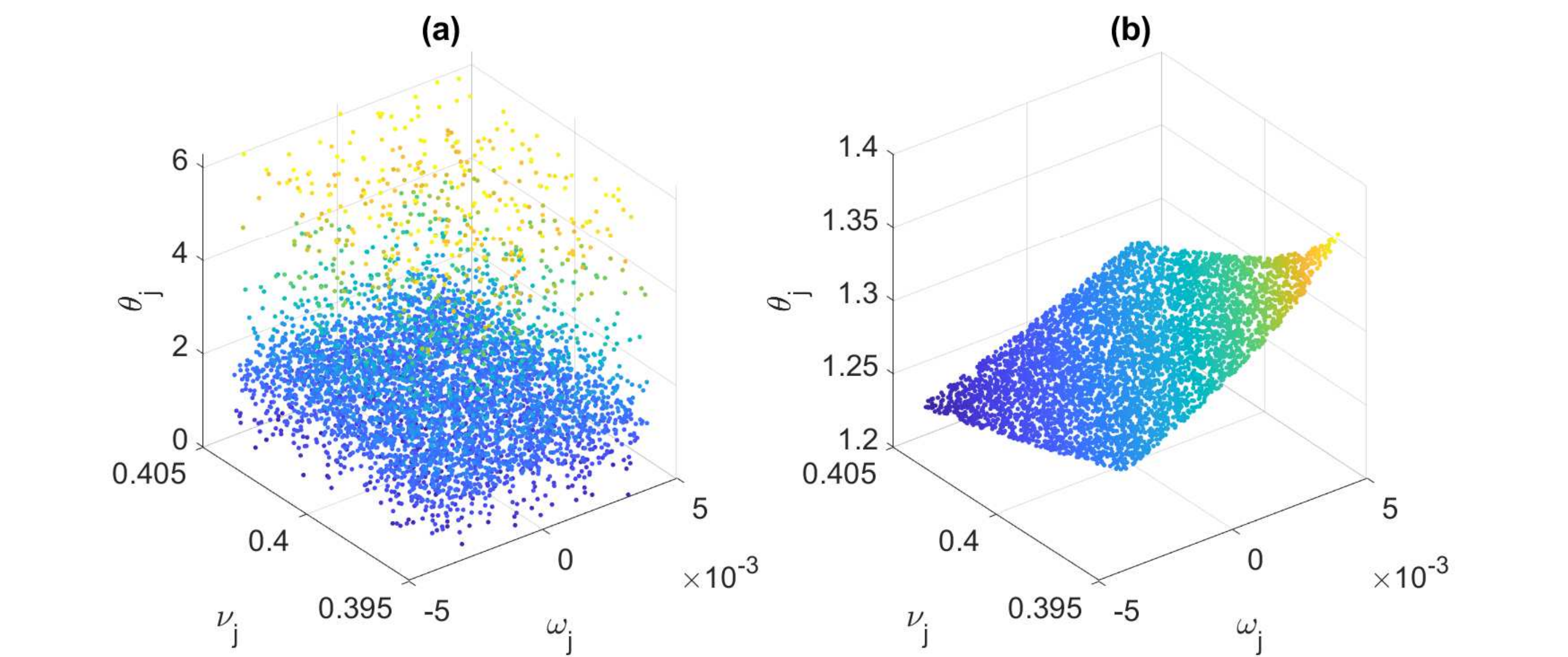}
\caption{Snapshot of a chimera state for a network of the form~\eqref{eq:dth1}-\eqref{eq:dth2}.
(a): population 1; (b): population 2. The colour shows the values of the $\theta_j$, with
different ranges in the two panels.
Parameters: $\mu=0.6,\nu_0=0.4,\Delta\nu=0.005,\omega_0=0,\Delta\omega=0.004,\beta=0.08,N=5000$.}
\label{fig:het2}
\end{center}
\end{figure}

While we have considered heterogeneity in only one parameter at a time, it is possible to consider
more than one. As an example, take a network of the form~\eqref{eq:dth1}-\eqref{eq:dth2}
where not only the $\omega_j$ are taken from a uniform distribution on 
$[\omega_0-\Delta\omega,\omega_0+\Delta\omega]$ but the values of $\nu$ are taken from a 
uniform distribution on $[\nu_0-\Delta\nu,\nu_0+\Delta\nu]$, as considered in
Sec.~\ref{sec:varynu}. A snapshot of a stable chimera state for such a network is shown
in Fig.~\ref{fig:het2}. We see that population 1 is incoherent while population 2 is
locked.  In the continuum limit the state of oscillators in population 2 could be described
by a function $\phi(\omega,\nu,t)$ defined for 
$(\omega,\nu)\in[\omega_0-\Delta\omega,\omega_0+\Delta\omega]\times [\nu_0-\Delta\nu,\nu_0+\Delta\nu]$, while those in population 1 would be described by a complex-valued function 
$a(\omega,\nu,t)$ defined for the same range of $(\omega,\nu)$. Numerical studies of such
systems would be more involved than the study of networks with a single heterogeneous parameter.

Another possibility is to consider nontrivial connectivity within or between populations. 
For example, we could replace~\eqref{eq:dth1} by
\be
   \frac{d\theta_j}{dt}  = \omega+\frac{\mu}{N}\sum_{k=1}^N\sin{(\theta_k-\theta_j-\alpha)}+\frac{\nu}{\langle d\rangle}\sum_{k=1}^N A_{jk}\sin{(\theta_{N+k}-\theta_j-\alpha)} \label{eq:dthmod}
\ee
and similarly for~\eqref{eq:dth2} where $A$ is the connectivity matrix between populations:
$A_{jk}=1$ if oscillator $k$ in one population is connected to oscillator $j$ in the other
and $A_{jk}=0$ otherwise. (Connections are undirected, so $A$ is symmetric.) $\langle d\rangle$
the mean degree: $\langle d\rangle=\sum_{j,k}A_{jk}/N$. If the connections between networks
are made randomly and each oscillator is connected to sufficiently many others we can make the
approximation~\cite{koerm08}
\be
   \frac{1}{\langle d\rangle}\sum_{k=1}^N A_{jk}\sin{(\theta_{N+k}-\theta_j-\alpha)}
\approx \frac{d_j}{N\langle d\rangle}\sum_{k=1}^N \sin{(\theta_{N+k}-\theta_j-\alpha)}
\ee
where $d_j$ is the degree of oscillator $j$: $d_j\equiv \sum_{k=1}^N A_{jk}$. Thus having a range
of degrees has approximately the same effect on the dynamics
as having a range of $\nu$ values, as investigated in
Sec.~\ref{sec:varynu}. To compare with the results in Sec.~\ref{sec:varynu} we construct networks
using the configuration model~\cite{new03}
having degree distributions which are uniform on $[d_0-\Delta d,d_0+\Delta d]$ and set $\nu=0.4$.
The corresponding value of $\Delta\nu$ as shown in Fig.~\ref{fig:vardelnukur} is
$\Delta\nu=\nu(\Delta d/d_0)=0.4(\Delta d/d_0)$, but note that both $\Delta d$ and $d_0$
are integers.

\begin{figure}
\begin{center}
\includegraphics[width=15cm]{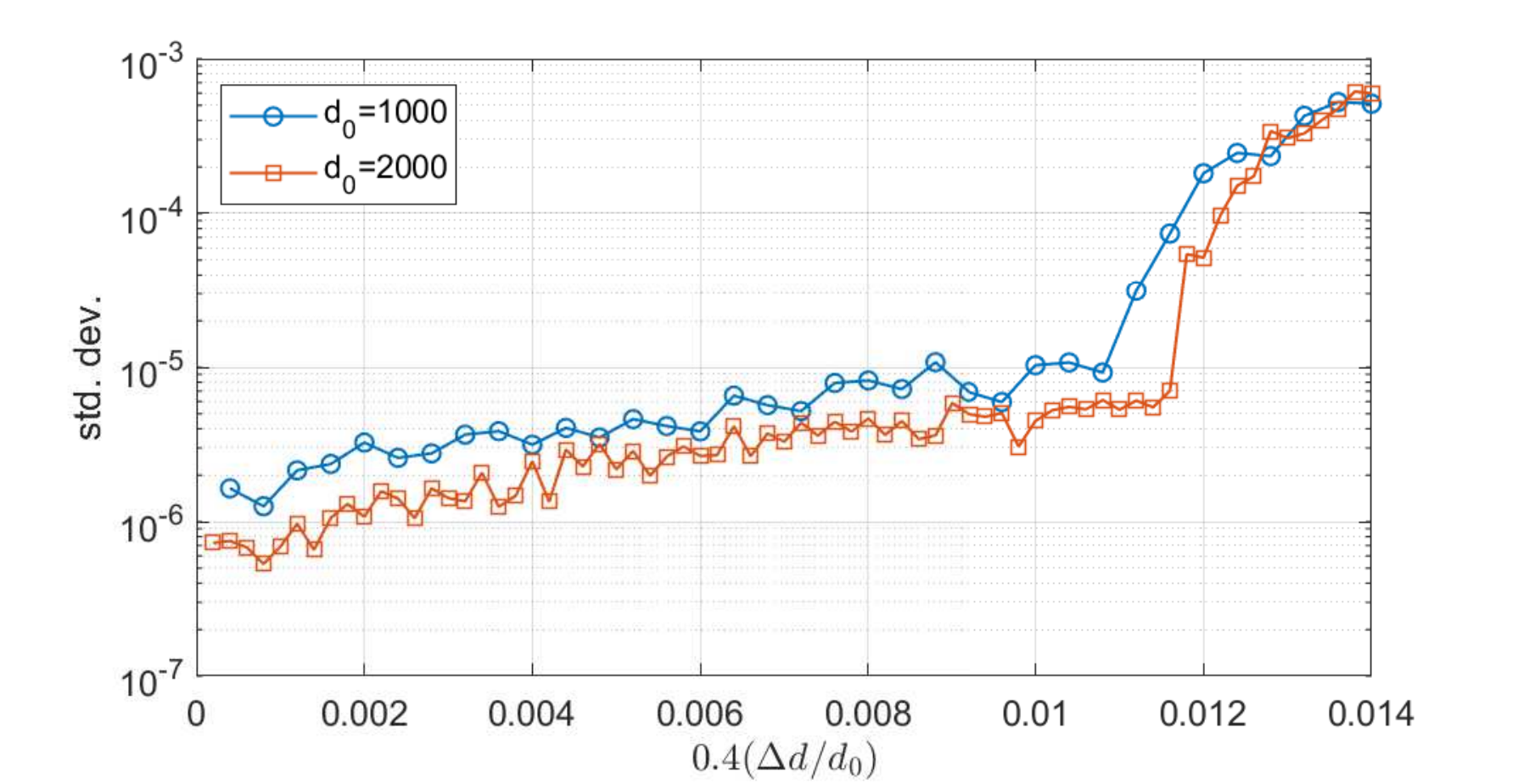}
\caption{Standard deviation of the order parameter~\eqref{eq:ord} as a function of the width
of the degree distribution, $\Delta d$.
Parameters: $\mu=0.6,\nu=0.4,\omega=0,\beta=0.08,N=5000$.}
\label{fig:degsc}
\end{center}
\end{figure}

We construct networks with $N=5000$ and $d_0$ either 1000 or 2000, and numerically 
solve equations of the form~\eqref{eq:dthmod} with initial conditions close to a chimera state.
We define the real order parameter for the synchronised (or largely synchronised) population
\be
   R=\left|\frac{1}{N}\sum_{j=1}^N e^{i\theta_j}\right| \label{eq:ord}
\ee
and (after transients) measure the standard deviation of $R$ over 300 time units. This 
standard deviation is
plotted on a log scale 
in Fig.~\ref{fig:degsc} as a function of $0.4(\Delta d/d_0)$ for the two different 
values of $d_0$. We see a rapid increase in the standard deviation indicating
the loss of full synchrony of the synchronised population when $0.4(\Delta d/d_0)$
is approximately 0.011-0.012, in very good agreement with Fig.~\ref{fig:vardelnukur}.
Note that several other authors have studied chimeras in a pair of subnetworks
with less than all-to-all connectivity~\cite{olmtor19,lairaj12,lisaa17}.

All of the results shown here have used a uniform distribution of heterogeneous parameters,
so to confirm the validity of these results we also considered a beta distribution with equal
shape parameters. For a beta distribution 
the distribution of a parameter $x$ with non-zero density on 
$[x_0-\Delta x,x_0+\Delta x]$ is proportional to
\be
   p(x)=\left(1-\frac{x-x_0}{\Delta x}\right)^\alpha\left(1+\frac{x-x_0}{\Delta x}\right)^\alpha
\ee
where $\alpha$ is the shape parameter. Integrals over $x$ are then approximated using
Gauss-Gegenbauer quadrature~\footnote{https://people.sc.fsu.edu/$\sim$ jburkardt/m\_src/gegenbauer\_rule/gegenbauer\_rule.html}.
We set $\alpha=2$ and used a discretisation of 10 points.
All of the results presented above for a uniform distribution were qualitatively reproduced
using this beta distribution (not shown).

%Figs 2,3,5,9,10,12 similar.

Although we considered networks with all-to-all connectivity within and between
populations, the governing equations were derived in the continuum limit. Thus the dynamics of large
but finite networks whose graphs have the same graphon (or graph limit)~\cite{lovsze06} as
the networks studied here, of the form
\be
   W(x,y)=\begin{cases} a, & (x,y)\in[0,1/2]\times[0,1/2] \mbox{ or }
 (x,y)\in[1/2,1]\times[1/2,1] \\
b, & \mbox{ otherwise}
\end{cases}
\ee
where $a$ and $b$ are constants with $b<a$ 
should also be described using the techniques presented here,
under the assumption that a parameter is uniformly distributed. Examples include 
Erd{\"os}-R{\'e}nyi networks, where the probability of connecting two oscillators within or between
populations is constant (but less than 1) and Paley graphs~\cite{chimed18}.

It might seem possible to generalise the techniques presented here to study chimeras on rings
of nonlocally coupled general oscillators~\cite{omezak15,omeome13}. However, while the locked
oscillators would be described by ODEs and the asynchronous ones by PDEs, one would need to know
(or to automatically find) the boundaries between such groups of oscillators, in order to
determine whether an ODE or a PDE was needed to describe the dynamics 
at a particular position on the ring. Also, these boundary points would move as parameters
were varied.

%(Could rewrite nicer) Rewrite:
%
%Thus
%\be
%   0  = \omega_{N+j}-\Omega+\frac{\mu}{N}\sum_{k=1}^N\sin{(\phi_{k}-\phi_j-\alpha)}+\frac{\nu}{N}\sum_{k=1}^N\sin{(\theta_{k}-\phi_j-\alpha)}
%\ee
%\be
%   =\omega_{N+j}-\Omega+\mu\left[\cos{\phi_j}\hat{S}-\sin{\phi_j}\hat{C}\right]
%+\nu\left[\cos{\phi_j}S-\sin{\phi_j}C\right]
%\ee
%where
%\be
%   \hat{S}\equiv \frac{1}{N}\sum_{k=1}^N\sin{(\phi_k-\alpha)};  \hat{C}\equiv \frac{1}{N}\sum_{k=1}^N\cos{(\phi_k-\alpha)};  S\equiv \frac{1}{N}\sum_{k=1}^N\sin{(\theta_k-\alpha)};  C\equiv \frac{1}{N}\sum_{k=1}^N\cos{(\theta_k-\alpha)}
%\ee
%
%Rotate locked population so that $\hat{S}=0$. Thus
%\be
%  0=\omega_{N+j}-\Omega-\sin{\phi_j}(\mu\hat{C}+\nu C)
%+\nu\cos{\phi_j}S
%\ee
%Averaging we obtain
%\[
%   0=\omega_0-\Omega-
%\]

%3D oscillators?

%Continuous spectrum only for systems with rotational invariance?

%Plot log of order parameters to compare with OA?
%\bibliographystyle{unsrt}
%\bibliography{chim}% Produces the bibliography via BibTeX.

\end{document}